\documentclass{amsart}

\usepackage{amscd,amssymb}
\def\user{{\lower 0.7 ex \hbox{\char'176}}}

\def\Z{{\mathbb Z}}
\def\Q{{\mathbb Q}}
\def\R{{\mathbb R}}
\def\C{{\mathbb C}}
\def\F{{\mathbb F}}

\def\P{{\mathbb P}}
\def\V{{\mathbb V}}
\def\E{{\mathbb E}}                    
\def\DD{{\mathbb D}}
\def\XX{{\mathbb X}}
\def\MM{{\mathbb M}}
\def\GG{{\mathbb G}}
\def\OO{{\mathbb O}}

\def\A{{\mathcal A}}
\def\B{{\mathcal B}}

\def\J{{\mathcal J}}
\def\M{{\mathcal M}}
\def\O{{\mathcal O}}
\def\T{{\mathcal T}}
\def\U{{\mathcal U}}
\def\X{{\mathcal X}}

\def\cG{{\mathcal G}}

\def\cD{{\mathcal D}}
\def\cR{{\mathcal R}}
\def\cS{{\mathcal S}}
\def\cZ{{\mathcal Z}}
\def\cRG{{\mathcal R}{\mathcal G}}

\def\Mbar{\overline{\M}}

\def\kappabar{\overline{\kappa}}
\def\taubar{\overline{\tau}}
\def\rhotilde{\tilde{\rho}}
\def\XXhat{\widehat{{\mathbb X}}}
\def\MMhat{\widehat{{\mathbb M}}}
\def\OOhat{\widehat{{\mathbb O}}}
\def\GGhat{\widehat{{\mathbb G}}}

\def\g{{\mathfrak g}}
\def\h{{\mathfrak h}}
\def\k{{\mathfrak k}}
\def\p{{\mathfrak p}}
\def\s{{\mathfrak s}}
\def\t{{\mathfrak t}}
\def\u{{\mathfrak u}}
\def\gbar{\overline{\g}}
\def\ubar{\overline{\u}}

\def\sp{\s\p}                             

\def\dot{{\bullet}}
\def\blank{\phantom{x}}
\def\lie{{\rm lie}}
\def\ass{{\rm ass}}
\def\com{{\rm com}}
\def\tor{{\rm tor}}
\def\pr{{\rm pr}}
\def\arf{{\rm arf}}

\def\a{{\mathbf \gamma}}

\def\G{\Gamma}
\def\w{\omega}

\newcommand\Rt{\operatorname{Root}} 
\newcommand\id{\operatorname{id}}
\newcommand\spec{\operatorname{Spec}}

\newcommand\Diff{\operatorname{Diff}}
\newcommand\Hom{\operatorname{Hom}}
\newcommand\End{\operatorname{End}}

\newcommand\Aut{\operatorname{Aut}}
\newcommand\Der{\operatorname{Der}}
\newcommand\Out{\operatorname{Out}}

\newcommand\out{\operatorname{out}}
\newcommand\Jac{\operatorname{Jac}}
\newcommand\Pic{\operatorname{Pic}}
\newcommand\Gr{\operatorname{Gr}}
\newcommand\Ray{\operatorname{Ray}}
\newcommand\CH{\operatorname{CH}}
\newcommand\Lie{\operatorname{Lie}}
\newcommand\Ass{\operatorname{Ass}}
\newcommand\Com{\operatorname{Com}}
\newcommand\ori{\operatorname{or}}
\newcommand\Symp{\operatorname{Sp}}
\newcommand\spine{\operatorname{spine}}

\newcommand\inn{\operatorname{int}}
\newcommand\initial{\operatorname{in}}
\newcommand\term{\operatorname{term}}

\newtheorem{theorem}{Theorem}[section]

\newtheorem{proposition}[theorem]{Proposition}
\newtheorem{corollary}[theorem]{Corollary}
\newtheorem{conjecture}[theorem]{Conjecture}

\theoremstyle{definition}

\newtheorem{example}[theorem]{Example}

\theoremstyle{remark}

\newtheorem{question}[theorem]{Question}
\newtheorem{problem}[theorem]{Problem}

\begin{document}

\title{Mapping Class Groups and Moduli Spaces of Curves}

\author{Richard Hain}

\address{Department of Mathematics\\ Duke University\\
Durham, NC 27708-0320, USA}

\email{hain@math.duke.edu}

\thanks{First author supported in part by grants from the National
Science Foundation}

\author{Eduard Looijenga}

\address{Faculteit Wiskunde en Informatica,
University of Utrecht, Postal Box 80.010, NL-3508 TA Utrecht,
The Netherlands}

\email{looijeng@math.ruu.nl}


\maketitle

\tableofcontents

\section{Introduction}
\label{sec:intro}

It is classical that there is a very strong relation between
the topology of $\M_g$, the moduli space of smooth projective curves
of genus $g$, and the structure of the mapping class group $\G_g$,
the group of homotopy classes of orientation preserving diffeomorphisms
of a compact orientable surface of genus $g$. The geometry of $\M_g$,
the topology of $\M_g$, and the structure of $\G_g$ are all intimately
related. Until recently, the principal tools for studying these topics
were Teichm\"uller theory (complex analysis and hyperbolic geometry),
algebraic geometry, and geometric topology. Recently, a fourth
cornerstone has been added, and that is physics which
enters through the theories of quantum gravity and conformal field theory.
Already these new ideas have had a remarkable impact on the subject through the
ideas of Witten and the work of Kontsevich. In this article, we survey some
recent developments in the understanding of moduli spaces. Some of these
are classical (do not use physical ideas), while others are modern. One
message we would like to convey is that algebraic geometers, topologists,
and physicists who work on moduli spaces of curves may have a lot to learn
from each other.

Having said this, we should immediately point out that,
partly due to our own limitations, there are important developments
that we have not included in this survey. Our most notable omission is
the arithmetic aspect of the theory, much of which originates in
Grothendieck's fundamental works \cite{groth:marche}, \cite{groth:esq}.
We direct readers to the volume \cite{groth:dessins} and to the recent
papers of Ihara, Nakamura and Oda for other recent developments (see
Nakamura's survey \cite{nakamura:survey} for references). Other
topics we have
not covered include conformal field theory and recent work 
of Ivanov \cite{ivanov:rigid}
and Ivanov and McCarthy \cite{ivanov-mccarthy} on homomorphisms
from mapping class groups and arithmetic groups to mapping class groups.
Of particular importance is Ivanov's version of Margulis rigidity
for mapping class groups \cite{ivanov:rigid} which he obtains using some
recent fundamental work of Kaimanovich and Masur \cite{masur} on the
ergodic theory of Teichm\"uller space.

We shall denote the moduli space of $n$ pointed smooth projective
curves of genus $g$ by $\M_g^n$. Knudsen, Mumford and Deligne
constructed a
canonical compactification $\Mbar_g^n$ of it. It is the moduli space
of stable $n$ pointed projective curves of genus $g$. It is a projective
variety with only finite quotient singularities. Perhaps the most important
developments of the decade concern the Chow rings%
\footnote{All Chow rings and cohomology groups in this paper are with
$\Q$ coefficients except when explicit coefficients are used.}
of $\M_g^n$ and $\Mbar_g^n$.

The first Chern class of the relative cotangent bundle of the universal
curve associated to the $i$th point is a class $\taubar_i$ in
$\CH^1(\Mbar_g^n)$. One can consider monomials in the $\taubar_i$'s of
polynomial degree equal to the
dimension of some $\Mbar_g^n$. For such a monomial, one can take the
degree of the monomial as a  zero cycle on $\Mbar_g^n$ to obtain a
rational number. These can be assembled into a generating function. Witten
conjectured that this formal power series satisfies a system of
partial differential operators.  Kontsevich proved this using topological
arguments, and thereby provided inductive formulas for these intersection
numbers. These developments are surveyed in Section~\ref{sec:ribbon}.

For each positive integer $i$, Mumford defined a {\it tautological class}
$\kappabar _i$
in $\CH^i(\Mbar_g)$. The restrictions $\kappa _i$ of these classes to
$\CH^\dot(\M_g)$ generate a subalgebra of $\CH^\dot(\M_g)$ which is called
the {\it tautological algebra} of $\M_g$. Faber has conjectured that this
ring has the structure of the $(p,p)$ part of the cohomology ring of a
smooth complex projective variety of complex dimension $g-2$. That is, it
satisfies
Poincar\'e
duality and has the ``Hard Lefschetz Property'' with respect to $\kappa _1$.
Considerable evidence now exists for this conjecture, much of which is
presented in Section~\ref{sec:chow}.
Other developments on the Chow ring, such as explicit computations in
low genus, are also surveyed there.

In the early 80s, Harer proved that the cohomology in a given degree
of $\M_g$ is independent of the genus once the genus
is sufficiently large relative to the degree. These stable cohomology
groups form a graded commutative algebra which is known to be free.
The tautological classes $\kappa_i$ freely generate a polynomial algebra
inside the stable cohomology ring. Mumford and others have conjectured
that the stable cohomology of $\M_g$ is generated by the $\kappa_i$'s.
Some progress has been made towards this conjecture which we survey
throughout the paper. In Section~\ref{sec:agstability}
we consider the stabilization maps from an algebro-geometric point of view,
and in Section~\ref{sec:algebras} we survey Kontsevich's methods for
constructing classes in the cohomology of the $\M_g^n$.

We have also tried to advertise the fecund work of Dennis Johnson
on the Torelli groups. The Torelli group $T_g$ is the subgroup of
the mapping class group $\G_g$ consisting of those diffeomorphism
classes that act trivially on the homology of the reference
surface. This mysterious group, in some sense, measures the
difference between curves and abelian varieties and appears to play
a subtle role in the geometry of $\M_g$. Johnson proved
that $T_g$ is finitely generated when $g\ge 3$ and computed its
first integral homology group. These computations have
direct geometric applications, especially when combined with M.~Saito's
work in Hodge theory --- for example, they restrict the
normal functions defined over $\M_g$ and its standard level covers.
{}From this, one can give a computation of the Picard group of the
generic curve with a level $l$ structure. Johnson's work and its
applications is surveyed in Section~\ref{sec:torelli}.

Since $\G_g$ is the orbifold fundamental group of $\M_g$, an algebraic
variety, one should be able to apply Hodge theory and Galois theory to study
its structure. In Section~\ref{sec:hodgemap} we survey recent work on
applications of Hodge theory to understanding the structure of the Torelli
groups, mainly via Malcev completion. In Section~\ref{sec:algebras}
we combine this Hodge theory with recent results of Kawazumi and Morita
to show that the cohomology of $\M_g$ constructed
by Kontsevich using graph cohomology are, after stabilization, polynomials
in the $\kappa_i$'s. Thus Hodge theory provides some evidence for Mumford's
conjecture that the stable cohomology of the
mapping class group is generated by the $\kappa_i$'s.

Some of the results we discuss have not yet appeared in the literature,
at least not in the form in which we present them. Rather than mention
all such results, we simply mention a few instances where we believe our
presentation to be novel: the correspondences in Section~\ref{subsec:correspondences},
the r\^ole of the fundamental normal function for orbifold fundamental groups in
Section~\ref{fundgroup}, Theorem~\ref{tautbound} and the contents of 
Section~\ref{subsec:relation}.
 
\medskip
\noindent{\it Notation and Conventions.}
All varieties will be defined over $\C$ unless explicitly stated
to the contrary. Unless explicit coefficients are used, all (co)homology
groups are with rational coefficients. We will often abbreviate {\it
mixed Hodge structure} by MHS. The sub- or superscript {\it pr}
on a (co)homology group will denote the primitive part in both the
context of the Hard Lefschetz Theorem and in the context of Hopf algebras.
\medskip

\noindent{\it Acknowledgements.} We would like to thank Carel Faber for his
comments on part of an earlier version of this paper and Shigeyuki Morita for explaining to  us some of his recent work. We also appreciate the useful
comments by a referee. 
We gratefully acknowledge support by the AMS that enabled us to attend this conference. 

\section{Mapping Class Groups}
\label{sec:groups}

Fix a compact connected oriented reference surface $S_g$ of genus $g$, and
a sequence of distinct points $(x_0,x_1,x_2,\dots )$ in $S_g$. Let us write
$S_g^n$ for the open surface $S-\{x_1,\dots ,x_n\}$ and
$\pi_g^n$ for its fundamental group $\pi_1(S_g^n ,x_0)$. This group admits a
presentation with generators $\alpha _{\pm 1},\dots ,\alpha_{\pm g},\beta_1,\dots
,\beta_n$ and relation
$$
(\alpha_1,\alpha_{-1})\cdots (\alpha_g,\alpha_{-g})=\beta_1\cdots\beta_n,
$$
where $(x,y)$ denotes the commutator of $x$ and $y$.\footnote{For $n=0$ the
righthand side is to be interpreted as the unit element.} The generators are
represented by loops that do not meet outside the base point;
$\beta_i$ is represented by a loop that follows an arc to a point close to
$x_i$, makes a simple loop around $x_i$, and returns to the base point along
the same arc.

Let $\Diff^+(S)^n_r$
denote the group of orientation preserving diffeomorphisms of $S$ that fix
the $x_i$ for $i=1,\dots ,n+r$, and are the identity on $T_{x_i}S$ for
$i=n+1,\dots ,n+r$. Although not really necessary at this stage, it is
convenient to assume that $2g-2+n+2r>0$. In other words, we do not consider
the cases where $(g,n,r)$ is $(0,0,0)$, $(0,1,0)$, $(0,0,1)$,
$(0,2,0)$ or $(1,0,0)$. We will keep this assumption throughout the paper.

The {\it mapping class group} $\Gamma_{g,r}^n$
is defined to be the group of connected components of this group:
$$
\Gamma_{g,r}^n = \pi_0 \Diff^+(S)^n_r.
$$
We omit the decorations $n$ and $r$ when they are zero.
The mapping class group  $\Gamma_g^n$ acts on $\pi_g^n$ by outer
automorphisms. A theorem that goes back to Baer (1928) and Nielsen (1927)
\cite{nielsen} identifies $\Gamma_g$, via this representation, with the
subgroup of $\Out (\pi_g)$ (of index two) that acts trivially on
$H_2(\pi_g)\cong H_2(S_g)$.

When $n\ge 1$ we can consider the diagonal action of $\Aut (\pi_g^n)$ on
$(\pi_g^n)^n$. Clearly, $\Out (\pi_g^n)$ acts on the set
of orbits of $\pi_g^n$ (which acts by inner automorphisms on each component) in
$(\pi_g^n)^n$.
Now $\Gamma_g^n$ can be identified with the group of outer automorphisms of
$\pi_g^n$ that preserve the image of $(\beta_1,\dots
,\beta_n)$ in $\pi_g^n\backslash (\pi_g^n)^n$.
If we choose $x_n$ as a base point, then a corresponding
assertion holds: $\Gamma_g^n$ can be identified with a subgroup of
$\Aut (\pi_1(S_g^{n-1} ,x_n))$ that is characterized in a similar way.
The evident homomorphism $\G_g^n\to \G_g^{n-1}$ is surjective and 
its kernel can be identified with $\pi_1(S_g^{n-1} ,x_n)$
(acting by inner automorphisms).  
Ivanov and McCarthy \cite{ivanov-mccarthy}
recently showed that the resulting exact sequence cannot be split.

\subsection{Generators and basic properties}\label{subsec:basic}
Although a lot is known about these groups they are still poorly
understood. Let us quickly review some of their basic
properties. Dehn proved in \cite{dehn} that the mapping class groups are
generated by the `twists' that are now named after him:
if $\alpha$ is a simple (unoriented) loop on
$S_g^{n+r}$, then parameterize a regular
neighborhood of $\alpha $ in $S_g^{n+r}$ by the cylinder
$[0,1]\times S^1$ (preserving orientations) and define an
automorphism of $S_g$ that on this neighborhood is given by
$(t,z)\mapsto (t ,e^{2\pi it}z)$ and is the identity elsewhere. The isotopy
class of this automorphism only depends on the isotopy class of
$\alpha $ and is called the {\it Dehn twist} along $\alpha$.
(Perhaps we should add that $\alpha$ is, in turn, already determined by
its free homotopy class, in other words, by the associated conjugacy
class in $\pi_g^n$.) The corresponding element of
$\G_{g,r}^n$ is the identity precisely when $\alpha$ bounds a disk in
$S_g-\{ x_{n+1},\dots ,x_{n+r}\}$ which meets $\{ x_1,\dots ,x_{n}\}$ in at
most one point. Several people have found a finite presentation for the
mapping class groups. One with few generators was given by Waynryb
\cite{waynryb}. From this presentation one sees that the mapping class
groups considered here are perfect when $g \ge 3$ (a result due to Powell
\cite{powell} in the undecorated case).

There is an obvious homomorphism
$\G_{g,n}\to \G_g^n$. It is easy to see that it is surjective and that the
kernel is generated by the Dehn twists around the points $x_1,\dots ,x_n$.
These Dehn twists generate a free abelian central subgroup of $\G_{g,n}$ of
rank $n$. Now recall that a central extension of a discrete group
$G$ by $\Z$ determines an extension class in $H^2(G;\Z )$; it has a
geometric interpretation
as a first Chern class. In the present case we have $n$ such
classes $\tau _i\in H^2(\G_g^n;\Z )$, $i=1,\dots ,n$.

Conversely, each subgroup of $H^2(G;\Z )$ determines a central
extension of $G$ by that subgroup. Harer proved that $H^2(\G _{g,r};\Z )$
is infinite cyclic if $g\ge 3$ \cite{harer:h2}, so that
there is a corresponding  central extension
$$
0\to\Z\to\widetilde\G _{g,r}\to\G _{g,r}\to 1.
$$
Since $H_1(\G_{g,r};\Z )$ vanishes, this central extension is perfect (and
universal).
A nice presentation of it was recently given by Gervais
\cite{gervais}. The (imperfect) central extension by $\frac{1}{12}\Z$
containing this extension appears in the theory of conformal blocks; it
has a simple geometric description which we will give in
Section~\ref{sec:moduli}.

\subsection{Stable cohomology}\label{subsec:stable}
The mapping class groups $\G_{g,r}^n$ turn up in a connected sum construction
that we describe next. It is convenient to do this in a somewhat abstract
setting.  Suppose we are given a closed, oriented (but not necessarily
connected) surface $S$, a finite subset $Y\subset S$, and a fixed point free
involution $\iota $ of $Y$. Assume that $\iota $ has been lifted to an
orientation reversing linear involution $\tilde\iota$ on the spaces of rays
$\Ray (TS|Y)$. The {\it real oriented blow up} $S_Y\to S$ is a surface
with boundary
canonically isomorphic to $\Ray (TS|Y)$. So $\tilde\iota$ defines an
orientation reversing involution of this boundary. Welding the boundary
components of $S_Y$ by means of this involution produces a closed surface
$S({\tilde\iota})$. Some care is needed to give it a differentiable structure
inducing the given one on $S_Y$. Although there is no unique way to do
this, all natural choices lie in the same isotopy class. If $S$ happens to
have a complex structure, then each choice of a real ray $L$ in
$T_pS\otimes_{\C} T_{\iota (p)}S$ determines a lift of $\iota$ over the
pair $\{p, \iota (p)\}$: if $l$ is a ray in $T_pS$, then $\tilde\iota (l)$
is determined uniquely by the condition $l\otimes_{\C}\tilde\iota (l)=L$.

If $S({\tilde\iota})$ is connected, then each finite subset
$X$ of $S-Y$ determines a natural homomorphism from the
mapping class group which is perhaps best denoted by $\G (S)_Y^X$
(a product of groups of the type $\G_{g,r}^n$) to the mapping class group
$\G (S(\tilde\iota))^X$. The image of this homomorphism is simply the
stabilizer of the simple loops indexed by $Y/\iota $ that are images of
boundary components of $S_Y$. Its kernel is a free abelian group whose
generators can be labeled by a system of representatives $R$ of $\iota$
orbits in $Y$. Indeed, for each element $y$ of $R$, take the composite of the
Dehn twist around $y$ and the inverse of the Dehn twist around $\iota (y)$.
These maps appear in the stability theorems and are at the root of the
recent operad theoretic approaches to the study of the cohomology of
mapping class groups.

\begin{theorem}[Stability theorem, Harer \cite{harer:stab}]
There exists a positive constant $c$ with the following property.
If $S(\tilde\iota )$ is connected and $S'$ is a connected component of $S$
and $X$ a finite subset of $S'\setminus Y$, then the homomorphism
$$
\G (S')_{Y\cap S'}^X\to \G (S(\tilde\iota ))^X
$$
induces an isomorphism on integral cohomology in degree $\le
c.\text{genus}(S')$.
\end{theorem}

The constant $c$ appearing in this theorem was $1/3$ in Harer's
original paper. It was later improved to $1/2$ by Ivanov in
\cite{ivanov:teichm}. Most recently, Harer \cite{harer:imp_stab} has showed
that we can take $c$ to be about $2/3$ and that this is the minimal
possible value. There is also a version for twisted coefficients, due to
Ivanov \cite{ivanov}.

Harer's theorem says essentially that the $k$th cohomology group of
$\G_{g,r}^n$ depends only on $n$, provided that $g$ is large
enough. These stable cohomology groups are the cohomology of a single group,
namely the group $\G_{\infty}^n$ of compactly supported mapping classes of a
surface $S_{\infty}$ of infinite genus (with one end, say) that fix a given
set of $n$ distinct points.

Among the homomorphisms defined above are maps
$\G_{g,1}\times \G_{g',1}\to \G_{g+g'}$. These stabilize and define
homomorphisms of
$\Q$ algebras
$$
\mu : H^{\bullet}(\G_{\infty})\to H^{\bullet}(\G_{\infty})\otimes
H^{\bullet}(\G_{\infty}).
$$
This defines a coproduct on $H^{\bullet}(\G_{\infty})$. Together
with the cup product, this gives $H^{\bullet}(\G_{\infty})$ the
structure of a connected graded-bicommutative Hopf algebra.
The classification of such Hopf algebras implies that
$H^{\bullet}(\G_{\infty})$ is free as a graded algebra and is
generated by its set of primitive elements
$$
H^{\bullet}_{\pr}(\G_{\infty}) := \{x\in H^+(\G_{\infty}):
\mu (x)=x\otimes 1+1\otimes x\}.
$$
For each $i>0$, Mumford \cite{mumford} and Morita \cite{morita:classes}
independently found a class $\kappa_i$ in $H^{2i}_{\pr}(\G_{\infty})$
(we shall recall the definition in Section~\ref{sec:agstability}) and
Miller \cite{miller} and Morita \cite{morita:classes} independently showed
that each $\kappa_i$ is nonzero. So the $\kappa_i$'s generate
a polynomial subalgebra of the stable cohomology.
Mumford conjectured that they span all of
$H^{\bullet}_{\pr}(\G_{\infty})$. This has been verified by Harer in a
series of papers \cite{harer:h2}, \cite{harer:h3}, \cite{harer:h4} in
degrees $\le 4$.\footnote{He also tells us that he has checked that there are
no stable primitive classes in degree 5.}

The first Chern class $\tau _i\in H^2(\G_g^n;\Z )$ stabilizes also and we may
think of it as an element of $H^{\bullet}(\G_{\infty}^n;\Z )$ ($i=1,\dots
,n$).
The forgetful map $\G_{\infty}^n\to \G_{\infty}$ gives
$H^{\bullet}(\G_{\infty}^n)$ the structure of a module over this Hopf
algebra. From the stability theorem one can deduce:

\begin{theorem}[Looijenga \cite{looijenga}]
The algebra
$H^{\bullet}(\G_{\infty}^n;\Z )$ is freely generated by the classes
$\tau _1,\dots ,\tau _n$ as a graded-commutative $H^{\bullet}(\G_{\infty};\Z
)$ algebra.
\end{theorem}

\section{Moduli Spaces}
\label{sec:moduli}

A conformal structure and an orientation on $S_g$ determine a complex
structure on $S_g$. The {\it Teichm\"uller space} $\X_{g,r}^n$ is
the space of conformal structures on $S_g$ (with some reasonable
topology) up to isotopies that
fix $\{ x_1,\dots ,x_{n+r}\}$ pointwise and act trivially on the tangent
spaces $T_{x_i}S$ for $i=n+1,\dots, n+r$.  It is, in a natural way, a complex
manifold of dimension $3g-3+n+2r$.
As a real manifold it is diffeomorphic to a cell. The group $\G_{g,r}^n$
acts naturally on it. This action is properly discontinuous and a subgroup of
finite index acts freely.

If $\G$ is any subgroup of $\G_{g,r}^n$ that
acts freely, then the orbit space $\G\backslash \X_{g,r}^n$ is a classifying
space for $\G$ and so its singular integral cohomology coincides with
$H^{\dot}(\G ;\Z )$. This is even true with twisted coefficients: if $V$
is a $\G$ module,
then the trivial sheaf over $\X_{g,r}^n$ with fiber $V$ comes with an obvious
(diagonal) action of $\G$. Passing to $\G$ orbits yields a locally constant
sheaf $\V$ on $\G\backslash \X_{g,r}^n$. The cohomology of this sheaf equals
$H^{\dot}(\G ;V)$. For an arbitrary subgroup $\G$ of $\G_{g,r}^n$ these
statements still hold as long as we take our coefficients to be $\Q$ vector
spaces (but $\V$ need no longer be locally constant). For $\G =\G_{g,r}^n$,
we denote the orbit space by $\M_{g,r}^n$.

The space $\M_{g,r}^n$ is, in a natural way, a normal analytic
space and the obvious forgetful maps such as $\M_{g,r}^n\to \M_g^n$ are
analytic. An interpretation as a coarse moduli space makes it possible
to lift this analytic structure to the algebraic category. To see this,
we first choose a nonzero
vector in each tangent space $T_{x_i}S_g$. Each triple $(C;x,v)$, where
$C$ is a connected nonsingular complex projective curve $C$ of genus $g$,
$x$ an
injective map $x:\{ 1,\dots ,n+r\} \to C$, and $v$ a nowhere zero section of
$TC$ over $\{n+1,\dots ,n+r\}$, determines an element of $\M_{g,r}^n$.
This point depends only on the isomorphism class of $(C,x,v)$ with
respect to the obvious notion of isomorphism. Since each conformal structure
on $S$ gives $S$ the structure of a
nonsingular complex projective curve, $\M_{g,r}^n$ can be identified can be
identified
with the space of isomorphism classes of such triples. From the work of
Knudsen, Mumford and Deligne, we know that $\M_g^n$ is, in a natural way, a
quasi-projective orbifold. Recall that they also constructed a projective
completion $\Mbar_g^n$ of $\M_g^n$, the {\it Deligne-Mumford completion}
\cite{del_mum}, that also admits the interpretation of a coarse moduli
space. Its points parameterize the connected stable $n$ pointed curves
$(C,x)$ of arithmetic genus $g$, where we now allow $C$ to have ordinary
double points, but still require $x$ to map to the smooth part of $C$ and
the automorphism group of $(C,x)$ to be finite. The {\it Deligne-Mumford
boundary} $\Mbar_g^n-\M_g^n$ is a normal crossing divisor in the orbifold
sense. There is a projective morphism $\Mbar_g^{n+1}\to\Mbar_g^n$, defined
by forgetting the last point. It comes with $n$ sections
$x_1,\dots ,x_n$.  The fibers of this morphism are stable $n$ pointed curves
(modulo finite automorphism groups) and the morphism can be regarded as the
universal stable $n$ pointed curve (in an
orbifold sense). Let $\omega$ denote the relative dualizing sheaf of
this morphism, considered as a line bundle in the orbifold sense.
We can then think of $\M_{g,r}^n$ as the set of $(v_{n+1},\dots ,v_{n+r})$
in the total space of $x_{n+1}^*\omega \oplus\cdots \oplus x_{n+r}^*\omega$
restricted to $\M_g^{n+r}$
that have each component nonzero. So $\M_{g,r}^n$ is also quasi-projective.

Each finite quotient group $G$ of $\G_g^n$ determines, in an obvious way,
 a Galois cover $\M_g^n[G]\to \M_g^n$. The Deligne-Mumford completion
$\Mbar_g^n[G]$ of this cover is, by definition, the normalization of
$\Mbar_g^n$ in $\M_g^n[G]$.

\begin{theorem}[Looijenga \cite{looijenga:cover}]
There exists a finite group $G$ such that
$\Mbar_g[G]$ is smooth with a normal crossing divisor as
Deligne-Mumford boundary.
\end{theorem}

This has been extended by De Jong and Pikaart \cite{jong}
to arbitrary characteristic,
and by Boggi and Pikaart (independently) to the $n$-pointed case.
(They show that it also can be arranged that each irreducible component of
the Deligne-Mumford boundary of $\Mbar_g^n[G]$ is smooth.)
This makes it relatively easy to define the Chow algebra of
$\Mbar_g^n$: if $\Mbar_g^n[G]$ is smooth, then define
$\CH^{\dot} (\Mbar _g^n)$ to be the $G$ invariant part of
$\CH^{\dot} (\Mbar _g^n[G])$ (we take algebraic cycles modulo rational
equivalence
with coefficients in $\Q$). It is easy to see that this is
independent of the choice of $G$.

The central extension of $\G _g$ by $\Z$ ($g\ge 3$) discussed in
Section~\ref{subsec:basic} takes the geometric form of a complex line
bundle over Teichm\"uller space with $\G _g$ action and hence yields an
orbifold line bundle over $\M_g$.
Its twelfth tensor power has a concrete description: it is the determinant
bundle of the direct image of the relative dualizing sheaf of
$\M_g^1\to\M_g$ (this is a rank $g$ vector bundle). The orbifold fundamental
group of the associated $\C^{\times}$ bundle is just the central extension
of $\G_g$ by $\frac{1}{12}\Z$ mentioned in Section~\ref{subsec:basic}.

Since $\G_g$ is perfect when $g\ge 3$, we have $H^1(\G _g)=0$. Ivanov has
asked the following question:

\begin{question}[Ivanov]
Is it true that $H^1(\G)$ vanishes for all finite index subgroups $\G$
of $\G_g$, at least when $g$ is sufficiently large?
\end{question}

This would imply that the Picard group of each finite unramified cover
of $\M_g$ (in the orbifold sense) is finitely generated. The answer
to Ivanov's questions is affirmative, for example, for subgroups of finite
index of $\G _g$, $g\ge 3$, that contain the Torelli group --- see
(\ref{van_h1}).

\section{Algebro-Geometric Stability}
\label{sec:agstability}
The Deligne-Mumford completion $\Mbar_g^n$ comes with a natural
stratification into orbifolds, with each stratum parameterizing stable
$n$ pointed curves of a fixed topological type $T$. Denote this stratum
by $\M (T)$. It has codimension equal to the number of singular points of $T$.
The
normalization of the topological type $T$ is an oriented closed surface $S$
that comes with $n$ distinct numbered points $X=\{ x_1,\dots ,x_n\}$ and a
finite subset $Y$ of $S-X$ with a fixed point free involution $\iota$,
so that $T$ is recovered by identifying the points of $Y$ according to
$\iota$.
These topological data define a  moduli space $\M (S)^{X\cup Y}$
of the same type (we hope that the notation is self-explanatory) and
there is a  natural morphism $\M (S)^{X\cup Y}\to \M(S/\iota )^X$ that is a
Galois cover of orbifolds. This morphism extends to a finite surjective
morphism from the Deligne-Mumford completion $\Mbar (S)^{X\cup Y}$ to the
closure of $\M (T)$ in $\Mbar_g^n$. The resulting morphism
$\Mbar (S)^{X\cup Y}\to\Mbar_g^n$ has
only self-intersections of normal crossing type and so carries a normal
bundle in the orbifold sense. This normal bundle is a direct sum of line
bundles with one summand for each $\iota$ orbit $\{ p,p'\}$, namely
$p^*\omega^{-1}\otimes p'{}^*\omega^{-1}$. (To see this, notice that
the restriction of the  universal curve to $\Mbar (T)$ has a quadratic
singularity along the locus defined by the pair  $\{ p,p'\}$. Associating
to a local defining equation its hessian determines a natural isomorphism
between $p^*\omega^{-1}\otimes p'{}^*\omega^{-1}$ and the normal bundle
of a divisor in the Deligne-Mumford boundary passing through $\Mbar(T)$.)
We now see before us an algebro-geometric
incarnation of the map  that appears in the stability theorem:
the set of normal vectors that point towards the interior $\M_g^n$
is the  restriction to $\M(S)$ of the total space of the direct sum of
$\C^{\times}$ bundles in this normal bundle. So $\M (S)_Y^X$ maps to the
latter space, and although we do not have a morphism $\M (S)_Y^X\to \M_g^n$,
the map on cohomology behaves as if there were. In particular, the map
$H^{\dot}(\M_g^n)\to H^{\dot}(\M (S)_Y^X)$  is a MHS morphism. So the stability
theorem implies:

\begin{theorem}
[Algebro-geometric stability]
Suppose that the finite set $X$ is contained in a connected component $S'$ of
$S$ of genus $g'$, so that $\M (S')^X_{Y\cap S'}$ appears as
a factor of $\M(S)_Y^X$.  Choose points in the remaining factors
so that we have an inclusion
of $\M (S')^X_{Y\cap S'}$ in $\M (S)_Y^X$. Then for
$k\le cg'$ the composite map
$$
H^k(\M_g^n)\to H^k(\M (S)_Y^X)\to H^k(\M (S')^X_{Y\cap S'})
$$
is an isomorphism and so is the map
$$
H^k(\M_{g'}^n)\to H^k(\M (S')^X_{Y\cap S'})
$$
induced by the forgetful morphism $\M (S')^X_{Y\cap S'}\to \M (S')^X\cong \M
_{g'}^n$.
These maps are also MHS morphisms.
\end{theorem}

So the stable rational cohomology $H^{\dot}(\G_{\infty}^n)$
comes with a natural MHS.
A geometric consequence of this result is that each stable rational
cohomology class of $\M_g^n$ (that is, a class whose degree is in
the stability
range) extends across the open part of the blow up of $\Mbar (T)$
parameterizing
the normal directions pointing towards the interior. Pikaart showed that
 these partial extensions can be made to come from a single extension
to $\Mbar_g^n$, at least if $g$ is large compared with $k$.
But then it is not hard to show that if this is possible for large $g$,
then it is possible in the stable range and so the conclusion is:

\begin{theorem}[Pikaart \cite{pikaart}]\label{pikaart_purity}
The restriction map $H^k(\Mbar_g^n)\to H^k(\M_g^n)$ is surjective
in the stable range. Consequently, the MHS on
$H^k(\G_{\infty}^n)$ is pure of weight $k$.
\end{theorem}

Mumford's Conjecture, if known, would imply this result, and so Pikaart's
Theorem is evidence for the truth of this conjecture.

We illustrate this theorem with the known stable classes.
We have seen in the previous section that $\Mbar_g^n$ comes with $n$
orbifold line bundles $x_i^*\omega$, $i=1,\dots ,n$. Let $\taubar_{n,i}$ denote
the first Chern class of this line bundle, regarded as an element of
$\CH^1(\Mbar_g^n)$. The restriction of this class to
$\CH^1(\M_g^n)$ is a pull-back of the restriction of $\taubar _{n-1,i}$ to
$\M_g^{n-1}$ (when $n\ge 1$) and so we denote that restriction simply by
$\tau_i$. The underlying cohomology class of $\tau_i$ in $H^{2i}(\M_g^n)\cong
H^{2i}(\G_g^n)$ is what we denoted earlier by that symbol, in particular, it is stable.

For the definition of the tautological classes of $\M _g^n$, we shall not use
Mumford's original definition, but a modification proposed by Arbarello-Cornalba. This might begin with the observation that 
the ``functor'' which associates to an $(n+1)$-pointed stable 
genus $g$ curve $(C;x_1,\dots ,x_n,x)$ the cotangent space $T_x^*C$
defines an orbifold line bundle over $\Mbar_g^{n+1}$. It is not quite the
same as the relative dualizing sheaf $\omega$ of the forgetful map
$\Mbar_g^{n+1}\to \Mbar_g^n$: a little computation shows that it is
in fact $\omega (\sum_{i=1}^n (x_i))$. This is perhaps a more natural
bundle to consider than $\omega$. In any case, we denote the direct image of
$c_1(\omega (\sum_{i=1}^n (x_i))^{i+1}\in\CH ^{i+1}(\Mbar_g^{n+1})$ under the
projection $\Mbar_g^{n+1}\to \Mbar_g^n$ by $\kappabar_{n,i}\in \CH^i(\M_g^n)$
and its restriction to $\M_g^n$ by $\kappa_{n,i}$. The cohomology class underlying  $\kappa_{n,i}$ can be regarded as an element of
$H^{2i}(\G^n_g)$ (of Hodge bidegree $(i,i)$). These cohomology classes
stabilize and, for $n=0$, they define the nonzero primitive elements of degree
$2i$ alluded to in \ref{subsec:stable}. 

We regard (for $k=0,1,\dots ,n$) $\CH ^{\dot}(\Mbar _g^n)$ as a
$\CH ^{\dot}(\Mbar_g^k)$-algebra via the obvious forgetful morphism, and
view the classes $\kappabar_{k,i}$ as elements of $\CH ^{\dot}(\Mbar_g^n)$
when appropriate. The class $\kappabar_{n,i}$ is then not equal to 
$\kappabar_{n-1,i}$, but according to formula (1.10) of \cite{arb_cor}
we have:
$$
\kappabar_{n,i}=\kappabar_{n-1,i} +(\taubar_{n,n})^i.
$$
As Arbarello-Cornalba explain, the classes $\kappabar_{n,i}$ possess
a nice property not enjoyed by Mumford's classes. First recall that
every stratum of $\Mbar_g^n$ is the image of a finite map
$\Mbar (S)^{X\cup Y}\to\Mbar_g^n$ and that $\Mbar (S)^{X\cup Y}$ is a
product of varieties of the type $\Mbar_{g_\alpha}^{n_\alpha}$. The
pull-back of $\kappabar_{n,i}$ along this map is the
sum  of the classes $\kappabar_{n_\alpha ,i}$ (pulled back along the
projection $\Mbar (S)^{X\cup Y} \to \Mbar_{g_\alpha }^{n_\alpha}$).
Carel Faber pointed out to us that a similar property is enjoyed by 
the divisor class of the Deligne-Mumford boundary, but we know of no other
examples. Since this behaviour is reminiscent of that of a primitive
element in a Hopf algebra under the coproduct, we ask:

\begin{question}
What other collections $\left\{\mu_{g,n}\in \CH ^k(\M_g^n)\right\}_{g,n}$
have this property?
\end{question}

\subsection{Correspondences between moduli
spaces}\label{subsec:correspondences}
There is an altogether different way to relate the cohomology of the
moduli spaces $\M_g^n$ for different values of $g$. This involves certain
Hecke type correspondences. For simplicity we shall restrict ourselves
to the undecorated case $n=0$. We return
to the reference surface $S_g$ and suppose that we are given a subgroup
$\pi$ of $\pi_g$ of finite index $d$, say.
(For what follows only its conjugacy class will matter.)
This subgroup determines an unramified
finite covering $\tilde S\to S$ of closed oriented surfaces. The genus
$\tilde g$ of $\tilde S$ is then equal to $d(g-1)+1$.
Consider the group of pairs $(\tilde h,h)\in \Diff^+(\tilde S)\times
\Diff^+(S)$ such that $\tilde h$ is a lift of $h$. Let $\G_g(\pi )$ be
its group of connected components. The projection $\G_g(\pi )\to \G_g$
has as kernel the group of covering transformations of $\tilde S \to S$
(so is finite) and its image consists of the outer automorphisms of $\pi_g$
that come from an automorphism which preserves the subgroup $\pi$
(so is of finite index, $e$, say). There is a corresponding finite covering
of moduli spaces $p_1:\M_g(\pi)\to \M_g$, where $\M_g(\pi )$ is simply the
coarse moduli space of finite unramified coverings of nonsingular complex
projective curves $\tilde C\to C$ topologically equivalent to $\tilde S\to S$.
There is also a finite map $p_2:\M_g(\pi)\to \M_{\tilde g}$. Together they
define a one-to-finite correspondence $p_2p_1^{-1}$ from $\M_g$ to
$\M_{\tilde g}$. This extends over
the Deligne-Mumford compactifications: if $p_1:\Mbar_g(\pi )\to \Mbar_g$
denotes the normalization of $\Mbar_g$ in $\M_g(\pi)$, then $p_2$ extends
to a finite morphism $p_2 :\Mbar_g(\pi)\to\Mbar_{\tilde g}$. We have
an induced map
$$
T_{\pi}:=e^{-1}p_{1*}p_2^*: \CH^{\dot}(\Mbar_{\tilde g})\to
\CH^{\dot}(\Mbar_g)
$$
and likewise on cohomology. A computation shows that any monomial in the
tautological classes is an ``eigen class'' for such correspondences:

\begin{proposition}
The map $T_{\pi}$ sends $\kappabar_{i_1}\kappabar_{i_2}\cdots
\kappabar_{i_r}$ to
$d^r\kappabar_{i_1}\kappabar_{i_2}\cdots \kappabar_{i_r}$.
\end{proposition}

This proposition suggests the consideration, for given positive integers
$r$ and $s$, of sequences of classes $(x_g\in\CH^s(\Mbar_g))_{g\ge 2}$ of
fixed degree that have the property that $T_{\pi}(x_g)=d^r x_{(d-1)g+1}$
for each index $d$ subgroup $\pi$ of $\pi_g$.

\begin{question}
Is for such a system the image of $x_g$ in $H^{\dot}(\M_g)$ stable?
Is it in fact a polynomial of degree $r$ in primitive stable classes?
\end{question}

An affirmative answer would give us a notion of stability for the Chow groups
of the moduli spaces $\M_g$.

\section{Chow Algebras and the Tautological Classes}
\label{sec:chow}

We have already encountered some of the basic classes on $\Mbar_g^n$:
the first Chern classes $\taubar_i\in \CH^1(\Mbar_g^n)$
($i=1,\dots ,n$) and the tautological classes $\kappabar_i$ ($i=1,2,\dots $).
More such classes come from the boundary: if $\prod _i\Mbar _{g_i}^{n_i}\to
\Mbar_g$ is a Galois covering of a stratum of the boundary as in
Section~\ref{sec:agstability},
then we can add to these the push-forwards along this map of the exterior
products of the corresponding
classes on the factors. Let us call the subalgebra of
$\CH^{\dot}(\Mbar_g^n)$ generated by all these classes the
{\it tautological subalgebra} and denote it by $\cR^{\dot}(\Mbar_g^n)$.
The image of this algebra in  $\CH^{\dot}(\M_g^n)$ is denoted by
$\cR^{\dot}(\M_g^n)$; it is generated by the classes
$\kappa_{n,i}$ ($i=1,2,\dots $) and $\tau _i$ ($i=1,\dots ,n$). It is possible
that these classes generate the rational Chow ring of $\Mbar_g^n$ modulo
homological equivalence, but this is of course unknown.
In any case, these subalgebras are preserved under pull-back and
push-forward along the natural maps that we have met so far.

The first computations
were done by Mumford \cite{mumford} who found a presentation of
$\CH^{\dot}(\Mbar_2)$.
Subsequently Faber \cite{faber} calculated $\CH^{\dot}(\Mbar_2^1)$,
$\CH^{\dot}(\Mbar_3)$ and obtained partial results
on $\CH^{\dot}(\Mbar_4)$. In all these cases the tautological algebra
is the whole Chow algebra.

This is also the case for $\Mbar_0^n$, whose Chow algebra was computed
by Keel. This is a very remarkable algebra which appears in other
contexts. Because of this, we describe it explicitly. We first introduce
notation for the divisor classes on $\Mbar_0^n$. 
The boundary divisor $\Mbar_0^n - \M_0^n$ parameterizes
all singular stable $n$ pointed rational curves. Its components
correspond to the topological types of $n$ pointed stable rational
curves with exactly one singular point.
Such curves have exactly two irreducible components.
By collecting the points $x_i$ lying on the same
component, we obtain a partition $P$ of $\{1,\dots ,n\}$ into two subsets.
The stability property implies that both members of $P$ have at least
two elements. We denote the corresponding class in $\CH^1(\Mbar_0^n)$ by
$D(P)$.

\begin{theorem}[Keel \cite{keel}] The Chow algebra $\CH^\dot(\Mbar_0^n)$
coincides with $H^\dot(\Mbar_0^n)$ and, as a $\Q$ algebra, is generated
by the $D(P)$'s subject to the following relations:
\begin{enumerate}
\item[(i)] If $\{i,j,k\}$ are distinct integers in $\{1,\dots ,n\}$,
then the sum of the $D(P)$'s for which $P$ separates $i$ from
$\{ j,k\}$ is independent of $j$ and $k$ (and equals $\taubar_i$).

\item[(ii)] $D(P)\cdot D(P')=0$ if $P$ and $P'$ are independent in the sense
that the partition they generate has four nonempty members.
\end{enumerate}
\end{theorem}

The relations (ii) are geometrically obvious since the divisors
$D(P)$ and $D(P')$ do not meet if $P$ and $P'$ are independent.
The additive relations (i) are not difficult to see either:
if $C$ is a stable $n$ pointed rational curve, then a moment of thought 
shows that there is a unique morphism $z:C\to \P^1$ that is an isomorphism on 
one irreducible component, constant on the 
other irreducible components, and is such that $z(x_i)=1$, $z(x_j)=0$ and 
$z(x_k)=\infty$.
The differential $z^{-1}dz$ restricted to $x_i$ defines a section
of $x^*_i\omega $. The image of $z^{-1}dz$ in $T^*_{x_i}C$  vanishes
precisely when $z$ collapses the irreducible component containing $x_i$.

In \cite{manin:trees}, Manin derives a formula for the
Poincar\'e polynomial of $\Mbar_0^n$. Such a formula was independently
found by Getzler \cite{getzler} who also obtained the
$\mathcal{S}_n$ equivariant Poincar\'e polynomial of $H^\dot(\Mbar_0^n)$.
That is, he determined the character
of the $\mathcal{S}_n$ representations $H^k(\Mbar_0^n)$,
$k\ge 0$. Kaufmann \cite{k_m:product} recently gave a formula for
the intersection number of classes of strata of complementary dimension.

We now turn to the Chow and cohomology algebras of the moduli spaces
$\M_g$. First we list some results about the Chow algebras.

\begin{enumerate}
\item[]$\CH^{\dot}(\M_1^n)=\Q$ for $n=1,2$ (folklore)\par
\item[]$\CH^{\dot}(\M_2  )=\Q$ (folklore),\par
\item[]$\CH^{\dot}(\M_2^1)=\Q [\tau]/(\tau^2)$ (Mumford \cite{mumford})\par
\item[]$\CH^{\dot}(\M_3  )=\Q [\kappa_1]/(\kappa_1^2)$ (Faber
\cite{faber}),\par
\item[]$\CH^{\dot}(\M^1_3)=\Q [\kappa_1,\tau ]/(\kappa_1^2,
4\tau ^2-\tau\kappa_1)$ (Faber\cite{faber}),\par
\item[]$\CH^{\dot}(\M_4  )=\Q [\kappa_1]/(\kappa_1^3)$ (Faber
\cite{faber}),\par
\item[]$\CH^{\dot}(\M_5  )=\Q [\kappa_1]/(\kappa_1^4)$ (Izadi \cite{izadi}
combined with Faber \cite{faber:hyp}).
\end{enumerate}

The reason that such computations can be made is that, when $g$ and
$n$ are both small, the moduli space $\M_g^n$ has a concrete
description. For example, when $g=2$, each curve is hyperelliptic
and therefore given by configuration of $6$ points
on the projective line. In the case $g=3$ a nonhyperelliptic curve is
realized by its canonical system as a quartic curve in $\P^2$. The
double cover of the projective plane along this curve is a Del Pezzo
surface of degree $2$, i.e., is obtained by blowing up $7$ points in
the plane in general position. General curves of genus 4 and 5 can be
described as complete intersections of multidegrees $(2,3)$ (in $\P^3$)
and $(2,2,2)$ (in $\P^4$), respectively.

\subsection{The tautological algebra of $\M_g$ and Faber's Conjecture}
\label{subsec:tautalg}
On the basis of  numerous calculations, Faber, around 1993, made the
following conjecture.

\begin{conjecture}[Faber\cite{seminar}]
The tautological algebra $\cR^{\dot}(\M_g)$ is a graded Frobenius algebra
with socle in degree $g-2$. That is, $\dim \cR ^{g-2}(\M_g)=1$, and the
intersection product defines a
nondegenerate bilinear form $\cR ^i(\M_g)\times \cR ^{g-2-i}(\M_g)\to
\cR ^{g-2}(\M_g)$ $(i=0,\dots ,g-2)$. Moreover, $\kappa_1$ has the
Lefschetz property in $\cR^{\dot}(\M_g)$ in the
sense that  multiplication by $(\kappa_1)^{g-2-2i}$ maps
$\cR^i(\M_g)$ isomorphically onto $\cR^{g-2-i}(\M_g)$ for $0\le i\le
(g-2)/2$.
\end{conjecture}

Since the conjecture was made, evidence for it has been growing. For example:

\begin{theorem}[Looijenga \cite{looijenga:taut}] The algebra
$\cR^{\dot}(\M_g)$ is trivial in degree
$>g-2$ and $\cR^{g-2}(\M_g)$ is generated by the class of the
hyperelliptic locus  (a closed irreducible variety of codimension
$g-2$).
\end{theorem}

In particular $\kappa_1^{g-1}=0$. Since $\kappa_1$ is ample on $\M_g$,
we recover a theorem of Diaz \cite{diaz} which asserts that every
complete subvariety of $\M_g$ must be of $\dim \le g-2$.

Actually, in \cite{looijenga:taut} a stronger result is proven, which, among
other
things, implies that 
$\cR^k(\M_g^n)=0$ for $k>g-2+n$. An induction argument then shows that
$\cR^{3g-3+n}(\Mbar_g^n)$ is spanned by the classes of the zero dimensional
strata. But zero dimensional strata can be connected by one dimensional strata
and the one dimensional strata are all rational. This shows that 
$\cR^{3g-3+n}(\Mbar_g^n)\cong\Q$.

Faber recently proved that the tautological class $\kappa_{g-2}$ is nonzero.
To describe his result, we find it convenient to
introduce a compactly supported version
of the tautological algebra: let $\cR^{\dot}_c(\M_g^n)$ be defined as the
set of elements in $\cR^{\dot}(\Mbar_g^n)$ that restrict trivially to the
Deligne-Mumford boundary. This is a graded ideal in
$\cR^{\dot}(\Mbar_g^n)$ and the intersection product defines a map
$$
\cR ^{\dot}(\M_g^n)\times\cR ^{\dot}_c(\M_g^n)\to\cR _c^{\dot}(\M_g^n)
$$
that makes $\cR _c^{\dot}(\M_g^n)$ a $\cR ^{\dot}(\M_g^n)$-module.
Notice that every complete subvariety of $\M_g$ of
codimension $d$ whose class is in $\cR ^{\dot}(\Mbar_g^n)$
defines a nonzero element of $\cR^d_c(\M_g^n)$ (but it is
by no means clear that such elements span $\cR^{\dot}_c(\M_g^n)$).
A somewhat stronger form of the first part of Faber's Conjecture is:

\begin{conjecture}\label{strongfaber}
The intersection pairings
$$
\cR^k(\M_g)\times \cR^{3g-3-k}_c(\M_g)\to \cR^{3g-3}_c(\M_g)\cong\Q ,
\quad k=0,1,2,\dots
$$
are perfect (Poincar\'e duality) and $\cR ^{\dot}_c(\M_g)$
is a free $\cR^{\dot}(\M_g)$ module of rank one.
\end{conjecture}

Faber \cite{faber:hyp} finds a compactly supported class
$I_g$ in $\cR_c^{2g-1}(\M_g)$ with $\kappa_{g-2}\cdot I_g\neq 0$.
So $\cR^{\dot}_c(\M_g)$ should be the ideal generated by this element.
Faber verified his conjecture for genera $\le 15$ by writing down many
relations in $\cR ^{\dot}(\M _g)$ (this evidently gives an upper bound)
and using the nonvanishing of $\kappa_{g-2}$ (this gives a surprisingly strong
lower bound). 

A refined form of Conjecture~\ref{strongfaber} (which we shall not state here) also takes care of the Lefschetz property.

\begin{question}
Does the tautological ring of $\Mbar _g^n$ satisfy Poincar\'e duality?
Does it  have the Lefschetz property with respect to $\kappabar_{n,1}$?
(It is known that $\kappabar_{n,1}$ is ample \cite{cor}.)
\end{question}

\subsection{Cohomology of some moduli spaces}\label{subsec:somecohom}
As may be expected, even less is known about the cohomology algebras.
Here is an incomplete list of special results. In genus 0 we have that
the Chow algebra of $\Mbar _0^n$ maps isomorphically onto its rational
cohomology algebra. The cohomology of $\M_0^n$ is easily computed if we
start out from the observation that this space is the projective 
arrangement of type $A_{n-2}$. It then follows for instance, that its 
cohomology in degree $p$ is of type $(p,p)$. There are similar
descriptions of the moduli spaces of $n$-pointed hyperelliptic curves of genus
$g$ when $n=0,1,2$ that involve arrangements of type $A$ or $D$. These again
should enable us to determine their rational cohomology ring, but it seems that
this hasn't been done yet. In the same spirit arrangements of various types
(among them $E_6$ and $E_7$) were used in \cite{looijenga:mthree} to prove that
$$
H^{\dot}(\M_3  )=\CH^{\dot}(\M_3)+\Q u,
$$
where $u$ is a class of degree $6$ of Hodge bidegree $(6,6)$ and 
$$
H^{\dot}(\M^1_3)=\CH^{\dot}(\M^1_3)+\Q u +\Q u\tau + \Q u\kappa _1+\Q v,
$$
where $v$ is a class of degree $7$ and of Hodge bidegree $(6,6)$. 

\begin{question} 
The image of the tautological algebra 
in $H^{2p}(\M_g^n)$ consists of classes of type $(p,p)$. Are all such classes
of this form?
\end{question}

A version of the Hodge conjecture asserts that the rational classes in degree
$2p$ of type $(p,p)$ are in the image of the Chow algebra, so modulo this conjecture we are asking whether
every Chow class on $\M_g^n$ is homologically equivalent to a tautological class.

\section{The Ribbon Graph Picture}
\label{sec:ribbon}

Around 1981 Thurston, Mumford and Harer observed that partial completions
of the Teichm\"uller spaces $\X _g^n$ with $n>0$ possess two natural
$\G _g^n$ equivariant triangulations. One is based on the hyperbolic
geometry of $S_g^n$ (Thurston) and the other based on the singular
euclidean geometry of $S_g^n$ (Mumford, Harer).
The last approach was actually a direct, but very powerful application of
work that Jenkins and Strebel had done 10--20 years earlier.
It is this approach that we shall explain.

The basic notion is that of a {\it ribbon graph}. This is a finite 
graph\footnote{For us a graph is a cell complex of pure dimension one; 
its zero cells are called {\it vertices} and its one cells {\it edges}. 
So it has no isolated vertices.} $G$ together with a cyclic order on the 
set of oriented edges\footnote{An oriented edge of a graph
is an edge together with an orientation of it.} emanating from each vertex. As we
shall see, there is
a canonical construction of a surface that contains $G$ and of which
$G$ is a deformation retract. This construction should explain the
name. We first give a somewhat more abstract
characterization of ribbon graphs which is very useful in some applications.

Let $X(G)$ be the set of oriented edges of $G$. Let $\sigma _1$ be
the involution of $X(G)$ that reverses the orientation of each edge. The
set $X_1(G)$ of $\sigma _1$ orbits can be identified with the set of edges
of $G$.
The cyclic orderings define another permutation $\sigma _0$ of $X(G)$
as follows. Each oriented edge $e$ has an initial vertex $\initial (e)$ and a
terminal vertex $\term (e)$. Define $\sigma _0(e)$ to be the successor of
$e$ with respect to the given cyclic order on the set of oriented edges that
have $\initial (e)$ as their initial vertex. The set of orbits
$X_0(G)$ of $\sigma _0$ can be identified with the set of  vertices of $G$.
Put $\sigma _{\infty}:= (\sigma _1\sigma _0)^{-1}=\sigma _0^{-1}\sigma _1$.
Call an orbit of this permutation a {\it boundary cycle}. (Draw a
picture to see why.) The set of boundary cycles wil be denoted
$X_{\infty}(G)$. These data form a complete invariant of $G$, for
we can reverse the construction
and associate to a finite nonempty set $X$ endowed with a fixed point free
involution $\sigma _1$ and a permutation $\sigma _0$ of $X$,
a ribbon graph $G(X,\sigma _0,\sigma _1)$ whose oriented edges are indexed
by $X$ and such that $\sigma_0$ and $\sigma_1$ are the permutations defined
above.

For every oriented edge $e$ of $G$ we form the one point compactification
$\Delta _e$ of the half strip $e\times [0,\infty )$; this is just a
$2$-simplex, parameterized in an unusual way. We make identifications
along the boundaries of these simplices with the help of $\sigma _0$
and $\sigma _{\infty}$:
$e\times \{0\}$ is identified with $\sigma _1e\times \{0\}$ and
and $\{\term (e)\}\times [0,\infty )$ with
$\{\initial (\sigma _{\infty}e)\}\times [0,\infty )$ (in either case, the
identification map is essentially the identity). This is easily seen to be
a compact, triangulated surface $S(G)$ that contains $G$ as a
subcomplex. Its vertex set can be identified with the disjoint union of
the vertex set of $G$ (so $X_0(G)$) and $X_{\infty}(G)$. We call vertices
of the latter type {\it cusps}. Notice that $G$ is a deformation retract
of $S(G)-X_{\infty}(G)$ and that the surface is canonically oriented
if we insist that the cyclic orderings of the edges emanating from each
vertex are induced by the orientation.
Let us say that the ribbon graph $G$ is {\it $n$-pointed} if we are given a injection $y:\{ 1,\dots ,n\}\hookrightarrow X_{\infty}(G)\cup X_0(G)$ whose image contains
$X_{\infty}(G)$ and the vertices of valency $\le 2$.

Suppose that we are given a {\it metric} $l$ on $G$. That is, a function
that assigns to every (unoriented) edge of $X$ a positive real number.
Give $[0,\infty )$ the standard metric and every half strip
$e\times [0,\infty )$ the product metric. This defines (at least locally) a
metric on $S(G)$. This metric is euclidean except possibly at the vertices.
However, it is not difficult to show that the underlying conformal structure
extends across all the vertices of $S(G)$ so that we end up with a compact
Riemann surface $C(G,l)$. Notice that each cusp has a ``circumference'' ---
this is the length of the associated boundary cycle. It is clear that we get
the same complex structure if $l$ is replaced by a positive multiple of it
and so we may just as well assume that the total length of $G$ is $1$. With
this convention, the sum of the circumferences of the cusps is $2$. 
The work of Jenkins and Strebel shows that all compact Riemann
surfaces arise in this way:

\begin{theorem}[Strebel \cite{strebel}]
Let $(C;x:\{ 1,\dots,n\}\hookrightarrow C)$ be an $n$-pointed connected Riemann surface (so that the complement of the image of $x$ has negative Euler
characteristic as usual) and let $c_1,\dots,c_n$ be nonnegative real numbers, not all zero. Then there exists an $n$-pointed metrized ribbon graph $(G,y,l)$, with $y(i)$ a cusp of $G$ of circumference $c_i$ when $c_i>0$ and a vertex of $G$ otherwise, such that 
$(C(G,l),y)$ and $(C,x)$ are isomorphic as $n$-pointed Riemann surfaces.
Moreover, $(G,y,l)$ is unique up to the obvious notion of isomorphism.
\end{theorem}

The results of Strebel also include a continuity property:
a continuous variation
of the complex structure on $C$ corresponds to a continuous variation of
$(G,y,l)$ in a sense that we make precise. Denote by $\cRG _g^n$ the set
of isomorphism classes of $n$-pointed ribbon graphs $(G,y)$ that are marked in the sense that we are given an isotopy class of homeomorphisms 
$h: S_g\to S(G)$ with $h(x_i)=y(i)$, $i=1,\dots ,n$. On this set $\G _g^n$ acts, and it is easy to see that the number of orbits of markings is finite.

Suppose that $(G,y,[h])$ represents an element of $\cRG _g^n$. Denote the
geometric realization of the
abstract simplex on the set $X_1(G)$ by $\Delta (G)$. Notice that the
metrics $l$ on $G$ that give $G$ unit length are parameterized by the
interior of $\Delta (G)$.  The circumferences of the cusps add up to two, so
half the cicumferences are the barycentric coordinates of 
a simplicial projection $\lambda :\Delta (G)\to\Delta ^{n-1}$.
Let $s$ be an edge of $G$ that is not a loop and does not connect two
vertices in the image of $y$. Then
collapsing that edge yields a member $(G/s ,y/s, [h]/s)$ of $\cRG _g^n$.
We can regard $\Delta (G/s)$ as a face of $\Delta (G)$. Making these
identifications produces a simplicial complex which we will denote by
$\XXhat _g^n$. It comes with a simplicial map
$\lambda:\XXhat _g^n\to \Delta ^{n-1}$. We have a simplicial action of $\G _g^n$ on $\XXhat _g^n$ which preserves the fibers of $\lambda$. The union of relative interiors of
simplices of $\XXhat _g^n$ indexed by the elements of $\cRG _g^n$ is
an open subset $\XX _g^n$ of $\XXhat _g^n$. The results of Strebel can
be strengthened to:

\begin{proposition}[cf.\ \cite{looijenga:cell}]
The above construction defines a $\G _g^n$ equivariant homeomorphism of
$\XX _g^n$ onto $\X _g^n\times\Delta ^{n-1}$.
\end{proposition}

Now consider the quotient space
$$
\MMhat _g^n:=\G _g^n\backslash\XXhat _g^n.
$$
This is a finite simplicial orbicomplex that is equipped with a simplicial
map $\lambda :\MMhat _g^n\to\Delta ^{n-1}$. We regard this complex as a
compactification of its open subset $\MM _g^n:=
\G _g^n\backslash\XX _g^n$. According to the above theorem,
the latter is canonically homeomorphic with $\M _g^n\times\Delta ^{n-1}$.
This raises the question of how this compactification compares to that
of Deligne-Mumford. The answer is essentially due to Kontsevich:

\begin{theorem}[Kontsevich \cite{kontsevich:airy}, see also
\cite{looijenga:cell}]
The simplicial orbicomplex $\MMhat _g^n$ is a
quotient space of $\Mbar _g^n\times\Delta ^{n-1}$. Moreover,
the part of $\MMhat _g^n$ where $\lambda _i>0$ carries an
oriented piecewise linear circle bundle
(in the orbifold sense) whose pull-back to
$\Mbar _g^n\times \{\lambda\in\Delta ^{n-1} | \lambda _i>0\}$
is the oriented circle bundle coming from the standard line bundle
$\taubar_i$.
In particular, the part of $\MMhat _g^n$ lying over the interior of $\Delta
^{n-1}$
carries the tautological cohomology classes underlying $\taubar _{n,i}$, $i=1,\dots ,n$.
\end{theorem}

The defining equivalence relation on
$\Mbar _g^n\times\Delta ^{n-1}$ is a little subtle and we refer to
\cite{looijenga:cell} for details regarding both statement and 
proof.\footnote{The space used by Kontsevich is not quite
$\MMhat _g^n$, but basically the part lying over
the interior of $\Delta ^{n-1}$ times a half line. In this case the
circumference map $\lambda$ has image $(0,\infty)^n$.}

This compactification of $\M _g^n\times\Delta ^{n-1}$ plays a crucial r\^ole
in Kontsevich's proof of the Witten conjectures. There are however
earlier applications. These include Harer's stability theorem we met before,
the computation of the Euler characteristic of $\M _g^n$, and the proof
that $\G _g^n$ is a virtual duality group of dimension
$4g-4+n$. We shall not explain the relation with stability here,
but we will briefly touch on the other applications.

There is also a remarkable arithmetic aspect of ribbon graphs that
is presently under intense investigation, but which
we merely mention in passing. This
is the observation, made by Grothendieck in a research proposal
\cite{groth:esq}, that for a metrized ribbon graph $(G,l)$ all of whose
edges have equal length, the corresponding Riemann surface $C(G,l)$ is,
in a canonical way, a ramified covering of the Riemann sphere $\P ^1$ with
ramification
locus contained in $\{ 0,1,\infty\}$. The graph $G$ appears here as the
preimage of the interval $[0,1]$,
its vertex set  as the preimage of $0$ and the set of cusps as the preimage
of $\infty$. The preimage of $1$ consists of the midpoints of the edges.
At these points we have simple ramification. A covering of this type is
naturally an algebraic curve defined over some number field.
Conversely, every connected covering of
the Riemann sphere of this type arises in this manner. The absolute Galois
group of $\Q$ acts on the collection of isomorphism types of such coverings,
and thus also on each finite set $\cR\cG _g^n$. It is very difficult to come
to grips with this action. For more information we refer to the collection
\cite{groth:dessins} and to Grothendieck's manuscripts \cite{groth:marche}
and \cite{groth:esq}.
Since these metrized ribbon graphs represent the barycenters of the simplices
of $\MM _g^n$, one can also think of this as an action of the absolute  Galois
group on the simplices of $\MM _g^n$, but the significance of this is not
clear to us.

\subsection{Virtual duality and virtual Euler characteristic}
\label{subsec:virtual}
We first make some observations about simplicial complexes. Let $K$ be a
simplicial complex, $L$ a subcomplex. Set $U:=K-L$. Then $U$ admits a
canonical deformation retraction onto the union of the closed
simplices of the barycentric subdivision of $K$ that lie in $U$.
This is a subcomplex, called the {\it spine} of $U$, whose $k$-simplices
correspond to strictly increasing chains
$\sigma _0\subset\sigma _1\subset\cdots\subset\sigma _k$
of simplices of $K$ not in $L$.
Further, if $\G$ is a group of automorphisms of $K$ that preserves $L$,
and if
\begin{enumerate}
\item[(i)] $U$ is contractible,
\item[(ii)] a subgroup of $\G$ of finite index acts freely on $U$,
\end{enumerate}
then $\G\backslash U$ is a simplicial `orbicomplex' that is also a virtual
classifying space for $\G$.\footnote{This means there is a normal subgroup 
$\G_1\subset\G$ of finite index such that a $\G /\G_1$-cover of this space classifies $\G_1$.} It has $\G\backslash \spine (U)$ as deformation
retract, and so the dimension of this spine is an upper
bound for the virtual homological dimension of $\G$.

We apply this in the situation where $K$ is the preimage of the first vertex of
$\Delta ^{n-1}$ in $\XXhat _g^n$ under $\lambda$ and $U=K\cap \XX _g^n$. Notice
that $U\cong \X _g^n$. The simplices meeting $U$ are indexed by the elements of
$\cRG _g^n$ with a single boundary cycle. A simple calculation shows that 
when $g\ge 1$, the number of edges of such a graph is at most $6g-5+2n$ and at 
least $2g-1+n$. For $g=0$ these numbers are $2n-5$, resp.\ $n-2$.
So the spine of $U$ has dimension $\le 4g-4+n$, resp. $n-3$. One can verify
that this is, in fact, an equality. It follows that $U$ admits a subcomplex of 
this dimension as an equivariant deformation retract. Hence:

\begin{theorem}[Harer \cite{harer:virt}] If $n\ge 1$, then  for every level
structure, the moduli space $\M _g^n[G]$ contains a subcomplex of dimension
$4g-4+n$ (when $g>0$) or $n-3$ (when $g=0$) as a deformation retract. 
\end{theorem}

From this he deduces a similar result for the case when $n=0$: $\M _g[G]$
has the homotopy type of complex of dimension $4g-5$.

\begin{problem}
Is there a Lefschetz type of proof of this fact? For instance, the Lefschetz
property would follow if one can find an orbifold stratification of $\M_g^n$  with all strata affine subvarieties of codimension $\le g$ ($n\ge 1$) 
or $\le g-1$ ($n=0$). That
would also show that the cohomological dimension of $\M_g^n$ for quasicoherent sheaves is $\le g-1$ ($n\ge 1$) or $\le g-2$ ($n=0$).
\end{problem}

Let us return to the general situation considered earlier and
suppose, in addition, that
\begin{enumerate}
\item[(iii)] $\G\backslash K$ is a finite complex, and
\item[(iv)] $U$ is a simplicial manifold of dimension $d$, say.
\end{enumerate}
These conditions are satisfied in the case at hand.
It is then natural (and standard) to assign to each simplex of $K$
the weight that is the reciprocal of the order of its $\G$ stabilizer.
This weighting is constant on orbits. Wall's Euler characteristic of $\G$
is simply the usual alternating sum of the number of $\G$ orbits of
simplices not in $L$, except that each is counted with its weight.
Equivalently, it is the orbifold Euler characteristic of the quotient
$\G\backslash K$.

In the present case, a ribbon graph $G$ defining a member of $\cRG _g^n$
gives a contribution $|\Aut (G)|^{-1}(-1)^{|X_1(G)|}$ to the virtual Euler
characteristic. The computation of the resulting sum is a combinatorial
problem that was first solved by Harer and Zagier. Kontsevich
\cite{kontsevich:airy} later gave a shorter proof. The answer is:

\begin{theorem}[Harer-Zagier \cite{harer_zagier}] The orbifold Euler
characteristic of $\M_g^n$ equals
$$
(-1)^{n-1}\,\frac{(2g+n-3)!}{(2g-2)!}\zeta (1-2g).
$$
Here $\zeta$ denotes the Riemann zeta function.
\end{theorem}

Harer and Zagier also find formulae for the actual Euler characteristics of
$\M _g^1$ and $\M _g$. These are often negative so that there must be lot
of cohomology in odd degrees.

For the discussion of virtual duality we go back to the general situation
and assume that beyond the four conditions already imposed we have:

\begin{enumerate}
\item[(v)] $L$ has the homotopy type of a bouquet of $r$-spheres.
\end{enumerate}

Then the theory of Bieri-Eckmann can be invoked in a virtual setting:
if we set $D:=\tilde H_r(L;\Z )$ and regard $D$ as a
$\G$ module in an obvious way, then $H_{d-r-1}(\G ;D)$ is of rank one
and for any $\G$-module $V$ with rational coefficients
the cap products
$$
\cap :H^k(\G ;V )\otimes H_{d-r-1}(\G ;D )\to
H_{d-r-1-k}(\G ;V\otimes D ),\quad k=0,1,2,\dots
$$
are isomorphisms. One calls $D$ the {\it Steinberg module} of $\G$.

Harer \cite{harer:virt} proves that in the present case hypothesis (v) is
satisfied: $L$ is a subcomplex of dimension $2g-3-n$, resp.\ $n-4$ which is
$(2g-4-n)$-connected, resp.\ $(n-5)$-connected when $g>0$, resp.\ $g=0$.

We shall call the corresponding orbifold local system $\DD$ over $\M_g^n$
the {\it Steinberg sheaf}. The homology group $H_{4g-4+n}(\M_g^n ;\DD )$
is of rank one. For every orbifold local system $\V$ of rational
vector spaces on $\M _g^n$, cap product with a generator of this homology group
defines isomorphisms
$$
H^k(\M_g^n ;\V )\stackrel{\sim}{\to}
H_{4g-4+n-k}(\M _g^n;\V\otimes\DD ),\quad k=0,1,2,\dots
$$
when $g>0$ and $n>0$ (and similar isomorphisms in the remaining cases).
In particular, taking $\V$ to be $\Q$, we see that $H_\dot(\M_g^n;\DD)$
has a canonical MHS. This suggests that $\DD$ has some Hodge theoretic
significance. Unfortunately it is not of finite rank, yet we wonder:

\begin{question}
Is the Steinberg sheaf motivic? In particular, does it have natural
completions that carry (compatible) Hodge and \'etale structures?
\end{question}

\subsection{Intersection numbers on the Deligne-Mumford completion}
\label{subsec:intersection}
The intersection numbers in question are those defined by monomials in
the $\taubar _i$'s. To be precise, define for
every such monomial $\taubar _1^{d_1}\dots \taubar _n^{d_n}$
(with all $d_i \ge 0$) the intersection number
$\int _{\Mbar_g^n} \taubar _1^{d_1}\dots \taubar _n^{d_n}$
where $g$ is chosen in such a way that this has a possibility of being
nonzero:
$3g-3+n=d_1+\cdots +d_n$. A physics interpretation suggests that we
should combine these numbers into the generating function
\begin{equation*}
\sum _{n=1}^{\infty} \frac{1}{n!}\sum _{g>1-{\frac{1}{2}} n}\,
\sum_{d_1 + \cdots +d_n = 3g-3+n} t_{d_1}\dots t_{d_n}
\int _{\Mbar_g^n} \taubar _1^{d_1}\cdots \taubar _n^{d_n}.
\end{equation*}
Now pass to a new set of variables $T_1,T_3,T_5,\dots$ by
setting
$$
t_i=1.3.5.\cdots(2i+1)T_{2i+1}.
$$
The resulting expansion $F(T_1,T_3,T_5,\dots )$ encodes all these
intersection numbers.

Witten \cite{witten}
conjectured two other characterizations of this function, both of which allow
computation of its coefficients. These were proved by Kontsevich in his
celebrated paper \cite{kontsevich:airy}. Perhaps the most useful
characterization is the one which says that $F$ is killed by a Lie algebra
 of differential operators isomorphic to the Lie algebra of polynomial vector
fields in one variable. This Lie algebra comes with a basis $(L_k)_{k\ge -1}$
corresponding to the vector fields $(z^k\partial /\partial z)_{k\ge -1}$ and
Witten verified the identities $L_k(F)=0$ for $k=-1,0$ within the realm of
algebraic geometry. However no such proof is known for $k\ge 1$. Kontsevich's
strategy is to represent the classes
$\taubar _1^{d_1}\cdots \taubar _n^{d_n}$ by piecewise differential forms
on the ribbon graph model that can actually be integrated.
This allows him to convert the intersection numbers into weighted sums over
ribbon graphs. This leads to a new characterization of the generating function
that is more manageable. Still a great deal of ingenuity is needed to
complete the proof of Witten's Conjecture.

\section{Torelli Groups and Moduli}
\label{sec:torelli}

In the early 80s, Dennis Johnson published a series of
pioneering papers \cite{johnson:fg,johnson:ker,johnson:h1} on the
Torelli groups. Although this work is in geometric topology, it has
several interesting applications to algebraic geometry. Here
we review some of his work.

First a remark on notation.
In the remainder of the paper we will write
$V_g$ for the symplectic vector space $H_1(S_g)$ and
$Sp_g(\Z)$ for the
group $\Aut (H_1(S_g,\Z),\langle \blank,\blank \rangle)$;
this does not really clash with standard notation, since a choice
of a symplectic basis of $H_1(S_g;\Z)$ identifies this with the standard
integral symplectic group of genus $g$. Likewise, $Sp_g$ will stand for
the algebraic $\Q$-group defined by the symplectic transformations of
$V_g$; so its group of $\Q$-points, $Sp_g(\Q )$, is just
the group of symplectic automorphisms of $V_g$.

The mapping class group $\Gamma_{g,r}^n$ acts on the homology of
the reference surface $S_g$. Since each of its elements preserves
the orientation of $S_g$, we have a homomorphism
\begin{equation}\label{homom}
\Gamma_{g,r}^n \to Sp_g(\Z).
\end{equation}
which is surjective. The {\it Torelli group} $T_{g,r}^n$ is defined to be
its kernel\footnote{Note that there is no general agreement on the
definition of $T_{g,r}^n$ when $r + n > 1$.} so that we have an extension
$$
1 \to T_{g,r}^n \to \Gamma_{g,r}^n \to Sp_g(\Z) \to 1.
$$
The homology groups of $T_{g,r}^n$ are therefore $Sp_g(\Z)$ modules.

The simplest kind of element of $T_{g,r}^n$ is a Dehn
twist along a simple loop in $S_g^{n+r}$ that separates
$S$ into two connected components. We call such a loop a {\it separating
simple loop}. Another type of element of $T_{g,r}^n$ is determined by a
{\it separating pair of simple loops}. This is a pair of two disjoint
nonisotopic loops $\alpha _1,\alpha _2$ on $S_g^{n+r}$
that together separate $S$ into two connected components. The Dehn
twist along
$\alpha _1$ composed with the inverse of the Dehn twist along $\alpha _2$
is in $T_{g,r}^n$. The first of Johnson's results is:

\begin{theorem}[Johnson \cite{johnson:fg,johnson:ker,johnson:h1}]
\label{johnson}
When $g\ge 3$, $T_{g,r}^n$ is generated by elements
associated to a finite number of separating simple loops and
a finite number of separating pairs of simple loops.
If $[S_g]\in \wedge ^2 H_1(S_g ;\Z)$ corresponds to the fundamental
class of $S_g$, then there are natural $Sp_g(\Z)$ equivariant surjective
homomorphisms
$$
\tau_g^1 : T_g^1 \to \wedge ^3 H_1(S_g;\Z ) \text{ and }
\tau_g : T_g \to \wedge^3 H_1(S_g;\Z )/([S_g]\wedge H_1(S_g;\Z )).
$$
In both cases, the kernel of $\tau$ is the subgroup generated by the
elements associated to simple separating loops. Finally, the kernels
of the induced homomorphisms
$$
H_1(T_g^1;\Z ) \to\wedge ^3 H_1(S_g ;\Z)\text{ and }
H_1(T_g^1;\Z ) \to\wedge ^3 H_1(S_g ;\Z )/([S_g]\wedge H_1(S_g ;\Z ))
$$
are both 2-torsion.
\end{theorem}

Johnson also finds an explicit description of this 2-torsion. We will
give it in a moment, but first we want to point out an algebro-geometric
consequence of this theorem. Let $\widetilde{\M}_g\subset\Mbar _g$ be the
complement of the irreducible divisor whose
generic point parametrizes irreducible singular stable curves, and let
$\widetilde{\M} _g^1$ be its preimage in $\Mbar _g^1$.

\begin{corollary}\label{cor:pioftildem}
When $g\ge 3$, the orbifold fundamental group of $\widetilde{\M}_g$
(resp.\ $\widetilde{\M}_g^1$) is isomorphic to an extension of $Sp_g(\Z )$
by $\wedge ^3H_1(S_g;\Z)/([S_g]\wedge H_1(S_g;\Z ))$ (resp.\
$\wedge ^3H_1(S_g;\Z )$).
\end{corollary}

Johnson's theorem shows that the $Sp_g(\Z)$ action on $H_1(T_{g,r}^n)$
is the restriction of a representation of the algebraic
group $Sp_g$. We shall see shortly the importance of this property.
Let $\lambda_1, \lambda_2,\dots ,\lambda_g$ be a
fundamental set of weights of $Sp_g$ so that $\lambda_j$ corresponds to
the $j$th fundamental representation of $Sp_g$. This last representation
can be realized as the natural $Sp_g$ action on the primitive part of
$\wedge^j V_g$.

The next result follows from Johnson's Theorem by standard arguments.

\begin{corollary}
For each $g\ge 3$, there is a natural $Sp_g(\Z)$ equivariant
isomorphism
$$
\tau_{g,r}^n : H^1(T_{g,r}^n) \stackrel{\sim}{\to}
V(\lambda_3) \oplus V(\lambda_1)^{\oplus(r+n)}.
$$
\end{corollary}

A theorem of Ragunathan \cite{ragunathan} implies that when $g\ge 2$,
the first cohomology of each finite index subgroup of $Sp_g(\Z)$ with
coefficients in a rational representation of $Sp_g(\Q)$ vanishes.
So Johnson's computation also gives:

\begin{corollary}\label{van_h1}
If $g\ge 3$, then every finite index subgroup of $\G_{g,r}^n$ that
contains $T_{g,r}^n$ has zero first Betti number.
\end{corollary}

The situation is very different when $g < 3$. The Torelli groups $T_1$
and $T_1^1$ are trivial, while Geoff Mess \cite{mess} proved that when
$g=2$, $T_2$ is a countably generated free group. He also computed
$H_1(T_2;\Z)$. It is the $Sp_2(\Z)$  module obtained by inducing the
trivial representation up to $Sp_2(\Z)$ from the stabilizer
$(\Z/2)\ltimes (SL_2(\Z)\times SL_2(\Z))$ of a decomposition of
$H_1(S_2;\Z)$ into two symplectic modules each of rank 2. (We shall
sketch a proof in the next subsection.) It is still unknown whether,
for any $g\ge 3$, $T_g$ is finitely presented.

\begin{problem}
Determine whether $T_g$ is finitely presented when $g$ is sufficiently large.
\end{problem}

Next, we describe Johnson's computation of the torsion in $H_1(T_g;\Z)$.
Denote the field of two elements by $\F_2$. Recall that an $\F_2$
quadratic form on $H_1(S_g;\F_2)$ associated to the mod two symplectic form
$\langle\blank ,\blank\rangle$ on $H_1(S_g;\F_2)$
is a function $\w : H_1(S_g;\F_2) \to \F_2$ satisfying
$$
\w(a+b) = \w(a) + \w(b) + \langle a,b\rangle .
$$
The difference between any two such is an element of $H^1(S_g;\F_2)$.
This makes the set $\Omega_g$ of such quadratic forms an affine space
over the $\F _2$ vector space $H^1(S_g;\F_2)$. Denote the algebra of
$\F _2$ valued functions on $\Omega_g$ by $S\,\Omega_g$. All such functions
are given by polynomials and so we have a filtration
$$
\F_2 = S_0\Omega_g \subset S_1\Omega_g \subset S_2\Omega_g \subset
\dots \subset S\,\Omega_g,
$$
where $S_d\Omega_g$ denotes the space of polynomial functions of degree
$\le d$.
Since $f=f^2$ for each $f\in S\,\Omega_g$, the associated graded algebra
is naturally isomorphic to the exterior algebra $\wedge^\dot H_1(S_g;\F_2)$.
The algebra  $S\,\Omega_g$ has as a distinguished element which is called
the {\it Arf invariant}, denoted here by $\arf$. If
$a_1,\dots, a_g, b_1,\dots,b_g$
is a symplectic basis of $H^1(S;\F_2)$, then $\arf$ is defined by
$$
\arf : \w \mapsto \sum_i \omega(a_i)\omega(b_i).
$$
It is an element of $S_2\Omega_g$, and its zero set $\Psi _g$ is an
affine quadric in $\Omega_g$. Let $S_d\Psi _g$ denote the image of
$S_d\Omega_g$ in the set of $\F_2$ valued functions on $\Psi_g$.

\begin{theorem}[Johnson \cite{johnson:h1}]
There are natural isomorphisms
\begin{gather*}
\sigma_{g,1} : H_1(T_{g,1};\F_2) \stackrel{\sim}{\to} S_3 \Omega_g,\quad
\sigma_g^1 : H_1(T_g^1;\F_2)\stackrel{\sim}{\to}
S_3\Omega_g/\F_2\,\arf ,\\
\sigma_g : H_1(T_g;\F_2)\stackrel{\sim}{\to} S_3\Psi_g
\end{gather*}
which are equivariant with respect to the $Sp_g(\F_2)$-action.
These induce natural isomorphisms
$$
H_1(T_{g,1};\Z)_\tor\cong S_2\Omega_g,\quad
H_1(T_g^1;\Z)_\tor \cong H_1(T_g;\Z)_\tor \cong S_2\Psi_g.
$$
Moreover, the natural isomorphisms
$$
\phi_{g,r}^n : H_1(T_{g,r}^n;\F_2)/H_1(T_{g,r}^n;\Z)_\tor
\stackrel{\sim}{\to}
\left[H_1(T_{g,r}^n;\Z)/\text{\rm torsion}\right]\otimes\F_2
$$
correspond, under the isomorphisms $\sigma_{g,r}^n$ and $\tau_{g,r}^n$,
to the obvious isomorphisms
\begin{gather*}
\phi_{g,1} : S_3\Omega_g/S_2\Omega_g\stackrel{\sim}{\to}
\wedge^3H_1(S_g;\F_2),\quad 
\phi_g^1 : S_3\Omega_g/(\F_2\,\arf + S_2\Omega_g) \stackrel{\sim}{\to}
\wedge^3H_1(S_g;\F_2),\\
\phi_g : S_3\Psi_g/S_2\Psi_g \stackrel{\sim}{\to}
\wedge^3H_1(S_g;\F_2)/([S_g]\wedge H_1(S_g;\F_2)).
\end{gather*}
\end{theorem}

The homomorphisms $\tau _g^1$ and $\tau _g$ admit direct conceptual
definitions that we will give later. Here we give a formula for the
image of the standard generators of $T_{g,1}$ in $\wedge^3H_1(S_g;\Z)$
and in $S_3\Omega _g$ under $\tau_{g,1}$ and $\sigma_{g,1}$, respectively.

Let $(\alpha _1,\alpha _2)$ be a separating pair of simple loops.
Let $t$ be the corresponding element of $T_{g,1}$ --- recall that this is
the product of the Dehn twist about $\alpha_1$ and the {\em inverse} of
the Dehn twist about
$\alpha_2$. The two loops decompose $S_g$ into two pieces $S'$ and $S''$, say,
where we suppose that $S'$ contains the point $x_1$.
We orient $\alpha _1$ and $\alpha _2$ as boundary components of $S''$.
The  resulting cycles are opposite in $H_1(S'';\Z)$: $[\alpha_2]=-[\alpha_1]$, 
and each spans the radical of the intersection pairing on
this group. So there is a well-defined element in
$\wedge ^2H_1(S'';\Z)/[\alpha_1]\wedge H_1(S'';\Z)$ representing the
intersection pairing on $H_1(S'';\Z)$. Its wedge with $[\alpha_1]$ can be
regarded as an element of $\wedge ^3H_1(S'';\Z)$. Since the inclusion
$S''\subset S_g$ induces an injection on first homology, we can also
view the latter as  an element of $\wedge ^3H_1(S_g;\Z)$. This is the
element $\tau _{g,1}(t)$;  it is clear that it only depends on the
image of $t$ in $T_g^1$.

Next we associate to $t$ a function
$\sigma_t:\Omega_g\to \F_2$ as follows. If $\omega\in\Omega_g$ takes the
value 1 on $[\alpha]$, then we put $\sigma_t(\omega )=0$; if it takes
the value 0 on $[\alpha]$, then the  restriction of $\omega$ to
$H_1(S'';\F_2)$ factors through a nondegenerate quadratic function
on $H_1(S'';\F_2)/\F_2[\alpha]$. Then $\sigma_{g,1}(t)(\omega)$ is its
Arf invariant. It can be shown that $\sigma_{g,1}(t)$ lies in $S_3\Omega$.

Now suppose that $t$ is the element of $T_{g,1}$ associated to a
separating simple loop $\alpha$. Denote the pieces $S'$ and $S''$ as before.
In this case, $\tau_{g,1}(t)$ is trivial and $\sigma_{g,1}(t)$ is the element 
of $S\Omega_g$ that assigns to $\omega$ the Arf invariant of its restriction to
$H_1(S'';\F_2)$. Notice that if $\alpha$ is a simple loop around $x_1$, then
$\sigma_{g,1}(t)$ is just the function $\arf$. (This explains why we mod out by
this function when passing from $T_{g,1}$ to $T_g^1$.) 

Without a base point there is no way of telling $S'$ and $S''$ apart. It is
because of this ambiguity that we have to restrict functions to $\Psi_g$ in
order to obtain a well defined function.

A diffeomorphism of $S_g$ onto a smooth projective curve $C$ determines
a natural isomorphism between $\Omega_g$ and the space of
theta characteristics of $C$ (i.e., square roots of the canonical bundle
$K_C$; see for instance Appendix B of \cite{acgh}). This suggests that
Johnson's computation should have
an algebro-geometric interpretation, if not interesting applications
to the geometry of curves.

\begin{problem}
Give an algebro-geometric construction of the epimorphism $T_g\to S_3\Omega$.
\end{problem}

Van Geemen has suggested such a construction (unpublished).

\subsection{Torelli space and period space}\label{subsec:torelli}
The group $T_g$ acts freely on $\X_g$. The quotient $\T_g$ is
therefore a complex manifold. It is called {\it Torelli space}. According
to the discussion at the beginning of Section~\ref{sec:moduli}, $\T_g$
is then a classifying space for $T_g$ so that there is a canonical isomorphism
$H_\dot(T_g;\Z ) \cong H_\dot(\T_g;\Z )$.
Torelli space has a moduli
interpretation; it is the moduli space of smooth projective curves
$C$ of genus $g$ together with a symplectic isomorphism
$$
\a : H_1(S_g;\Z) \to H_1(C;\Z).
$$
There are also decorated versions $\T_{g,r}^n$ of Torelli space. Their
points are points of $\M_{g,r}^n$ together with a symplectic
isomorphism $\a$ of $H_1(S;\Z)$ with the first homology of the curve
corresponding to the point of $\M_g$. It is clear that
the map $\T_{g,r}^n \to \M_{g,r}^n$ is Galois with Galois group
$Sp_g(\Z)$.

Denote the Siegel space associated to $V_g$ by $\h_g$. To be precise,
$\h_g$ is the set of pure Hodge structures on $V_g$
with Hodge numbers $(-1,0)$ and $(0,-1)$, polarized by the intersection form.
This is a contractible complex manifold of dimension $g(g+1)/2$ on
which the group $Sp_g(\R)$ acts properly and transitively.
We can also regard  $\h_g$ as the moduli space of pairs consisting of a
$g$ dimensional principally polarized abelian
variety $A$ plus a symplectic isomorphism
$$
\a : H_1(S;\Z) \to H_1(A;\Z).
$$
This interprets the $Sp_g(\Z)$ orbit space of $\h_g$ as the moduli space
of principally polarized abelian varieties of dimension $g$, $\A _g$. We
regard $\A _g$ as an orbifold with orbifold fundamental group $Sp_g(\Z)$,
although $Sp_g(\Z)$ does not act faithfully on $\h_g$. The kernel of this
action is $\{\pm 1\}$.

Assigning to a smooth projective curve the Hodge structure on
its first homology group defines
a map $\T_g \to \h_g$, the {\it period map} for $T_g$. It is an isomorphism
in genus 1, an open imbedding when $g=2$, and 2:1 with ramification along
the hyperelliptic locus when $g\ge 3$.%
\footnote{It is stated incorrectly in \cite{hain:normal} that $\T_2 \to
\h_2$ is an unramified 2:1 map onto its image.}
The reason for this is that
for all abelian varieties we have the equality
$$
[A;\a] = [A;-\a]
$$
of points of $\h_g$ as $-\id$ is an automorphism of each abelian variety.
On the other hand, we have the equality
$$
[C;\a] = [C;-\a]
$$
of points of $\T_g$ if and only if $C$ is hyperelliptic.

Mess's result (mentioned at the beginning of the section)
can now be deduced from this: $\T_2$ is the complement
in $\h_2$ of the locus of principally polarized abelian varieties
that are products of two elliptic curves. The locus of such reducible
abelian varieties is a countable disjoint union of copies of
$\h_1 \times \h_1$. The group $Sp_2(\Z)$ permutes them transitively,
and each is stabilized by a product of two copies of $SL_2(\Z)$ and
an involution that switches the two copies of the upper half plane.
Mess's result follows easily using the stratified Morse theory of
Goresky and MacPherson --- use distance from a generic point of
$\h_2$ as the Morse function. Since each component of $\h_2 -\T_2$ is
a totally geodesic divisor, the distance function has a unique critical
point (necessarily a minimum) on each stratum. It follows that $\T_2$
has the homotopy type of a wedge of circles, one for each component of
$\h_2 -\T_2$.

The period map gives, after passage to $Sp _g(\Z )$ orbit spaces, a morphism
$\M_g\to\A_g$, the period mapping for $\M _g$. This period mapping extends
 to the partial completion
$\widetilde{\M}_g$ of $\M_g$ and the resulting map
$\widetilde{\M}_g\to \A_g$ is proper.

Now assume $g\ge 3$ and
denote the image of the period map $\T_g \to \h_g$ by $\cS_g$. This
space is the quotient of $\T _g$ by the subgroup $\{\pm 1\}$ of $Sp_g(\Z)$.
Consequently
$$
H^\dot(\cS_g) \cong H^\dot(T_g)^{\{\pm 1\}}.
$$
Observe that $\cS_g$ is a locally closed analytic subvariety of
$\h_g$, but not closed. The $\{\pm 1\}$ cover $\T_g\to\cS_g$ extends as
a $\{\pm 1\}$ cover $\overline{\T}_g\to\overline{\cS}_g$ over the closure
of $\cS_g$ in $\h _g$, and the $\{\pm 1\}$ action on the total space
is the restriction of an $Sp_g(\Z)$ action. Both $\overline{\T}_g$ and
$\overline{\cS}_g$ are rather singular along the added locus
(which is of codimension $3$). If we pass to $Sp_g(\Z)$ orbit spaces, then
the natural map
$$
\widetilde{\M}_g\to Sp_g(\Z)\backslash \overline{\T}_g\cong
Sp_g(\Z)\backslash \overline{\cS}_g
$$
resolves these singularities in an orbifold sense. A resolution of a normal analytic variety always induces a surjection on fundamental groups and so it follows  from (\ref{johnson}) that the fundamental group of $\overline{\T}_g$ is abelian  and is $Sp_g(\Z)$ equivariantly
a quotient of $\wedge ^3 H_1(S_g ;\Z )/([S_g]\wedge H_1(S_g ;\Z ))$.

\begin{problem}
Understand the topology of $\cS_g$ and its closure $\overline{\cS}_g$ in
$\h_g$. In particular, how close is $\cS_g$ to being a finite complex? (Observe
that if it has a finite 2-skeleton, then $T_g$ is finitely presented.)
\end{problem}

Related, but formally independent of this problem, is the question of
whether the cohomology of $T_g$ stabilizes in a suitable sense:

\begin{question}
Is $H^k(T_g)$ expressible as an $Sp_g(\Z )$ module in a manner that is
independent of $g$ if $g$ is large enough? For example, from Johnson's
Theorem, we know that $H^1(T_g)$ is the third fundamental representation
of $Sp_g$ for all $g\ge 3$.
\end{question}

\subsection{The Johnson homomorphism}\label{subsec:johnson}
The proof of Johnson's Theorem is non-trivial and uses geometric
topology, but the homomorphism $\tau_g^1$ is easily described.

Since $T_g$ is torsion free, the projection $\T _g^1\to \T_g$ defines the
universal curve over $\T_g$. Denote the corresponding bundle of jacobians by
$\J_g \to \T_g$. Since the local system of first homology groups
associated to the universal curve is canonically framed, this jacobian
bundle $\J_g \to \T_g$ is analytically trivial as a bundle of Lie groups: we
have a natural trivializing projection  $p:\J _g \to \Jac S_g$, where
$\Jac S_g:= H_1(S_g; \R /\Z )$ is the ``jacobian'' of the reference surface.

The usual Abel-Jacobi map, which assigns to an ordered pair of points $(x,y)$
on
a smooth curve $C$ the divisor class of $(x)-(y)$, induces a morphism
$$
\T_g^1\times _{\T _g} \T_g^1\to \J _g.
$$
over $\T _g$. This provides a correspondence
$$
\begin{CD}
\T_g^1\times _{\T _g} \T_g^1 @>>> \J _g @>p>> \Jac S_g  \cr
@VV{pr_2}V \cr
\T_g^1 \cr
\end{CD}
$$
from $\T_g^1$ to $\Jac S_g $. It induces homomorphisms
$$
H_k(T_g^1) \cong H_k(\T_g^1) \to H_{k+2}(\Jac S_g ).
$$
The first of these is the Johnson homomorphism
$$
\tau_g^1 : H_1(T_g^1) \to H_3(\Jac S_g )
$$
for $T_g^1$.
Since $\T _g^1\to \T_g$ is a fibration of Eilenberg-MacLane spaces, we have
an exact sequence of fundamental groups:
$$
1 \to \pi_g \to T_g^1 \to T_g \to 1.
$$
This induces an exact sequence
$$
H_1(S_g;\Z ) \to H_1(T_g^1;\Z ) \to H_1(T_g;\Z ) \to 0
$$
on homology. Since $\Jac S_g $ is a topological group with torsion free
homology, its integral homology has a product --- the {\it Pontrjagin
product}. It is not difficult to check that the composite
$$
H_1(S_g;\Z ) \to H_1(T_g^1;\Z ) \to H_3(\Jac S_g;\Z )
$$
is the map given by Pontrjagin product with the class $[S_g]$.
It follows that there is a natural homomorphism
$$
H_1(T_g;\Z ) \to H_3(\Jac S_g;\Z )/\left([S_g ]\times H_1(S_g;\Z )\right).
$$
This is the Johnson homomorphism $\tau_g$ for $T_g$.

Johnson's Theorem, alone and in concert with Saito's theory of
Hodge modules, has several interesting applications to the
geometry of moduli spaces of curves as we shall see in subsequent
sections.

\subsection{Monodromy of roots of the canonical bundle}
In this subsection we assume that $g\ge 2$. Suppose that $C$ is a smooth
projective curve of genus $g$. Since its canonical bundle $K_C$ is of degree
$2g-2$ and since $\Pic^0 C$ is a divisible group, $K_C$ has $n$th roots
whenever $n$ divides $2g-2$. Any two such $n$th roots will differ by an
$n$ torsion point of $\Pic^0 C$. Because of this, $n$th roots of $K_C$
are rigid under deformation. It follows that they form a locally constant
sheaf (in the orbifold sense) $\Rt^n$ over $\M_g$. The fiber over
$C$, denoted $\Rt^n C$, is a principal homogenous space over
$H_1(C;\Z/n)$, the group of $n$ torsion points of $\Pic^0 C$.

Choose a conformal structure on $S_g$. Denote the corresponding
algebraic curve by $C$.
Sipe \cite{sipe} determined the monodromy representation
$$
\rho^n : \G_g \to \Aut \Rt^n(C).
$$
of this sheaf. Before giving it, we make some remarks. Since the Torelli
group acts trivially on the $n$ torsion of $\Pic^0 C$, it follows that
the restriction of $\rho^n$ to $T_g$ factors through a representation
$T_g \to H_1(S;\Z/(2g-2)) \to H_1(S;\Z/n)$. However, the action of $\G_g$ on
the set $\Rt^2 C$ of square roots of $K_C$ (the set of theta characteristics
of $C$) factorizes through $Sp_g(\Z)$ also (even through $Sp_g(\F_2)$) --- this
is because there is a canonical correspondence between square roots of $K_C$
and $\F_2$ quadratic forms on $H_1(S;\F_2)$ associated to the intersection form.
It follows that the image of the monodromy representation $\rho^n$ will
be contained in an extension of $Sp_g(\Z/n)$ by a subgroup of
$2\cdot H^1(S_g;\Z/n)$. In fact, it is all of this group.

\begin{theorem}[Sipe \cite{sipe}]\label{sipe}
The monodromy group of $\Rt^n$ is an extension of
$Sp_g(\Z/n)$ by the subgroup $2\cdot H^1(S_g;\Z/n)$ of $H_1(S_g;\Z/n)$.
\end{theorem}

The subgroup $2\cdot H^1(S_g;\Z/n)$ appears as a quotient of the Torelli group
$T_g$. In \cite{hain:normal} it is shown that the restriction of the monodromy
representation to $T_g$ is the composite of the Johnson homomorphism with a natural surjection
$$
\wedge^3H_1(S_g;\Z)/([S_g]\wedge H_1(S_g;\Z)\to H_1(S_g;\Z/(g-1)) \to 2\cdot H_1(S_g;\Z/n).
$$

\subsection{Picard groups of level covers}\label{subsec:picard}
Denote the moduli space of smooth projective genus $g$
curves with a level $l$ structure by $\M_g[l]$. This is convenient
shorthand for the notation $\M_g[Sp_g(\Z/l)]$ introduced in
Section~\ref{sec:moduli}. Denote the kernel of the reduction mod $l$ map
$$
Sp_g(\Z) \to Sp_g(\Z/l)
$$
by $Sp_g(\Z)[l]$, and its full inverse image in $\Gamma_g$ by
$\Gamma_g[l]$. Then $\M_g[l]$ is the quotient of Teichm\"uller
space $\X_g$ by $\Gamma_g[l]$. As in the case of $\M_g$, there is
a canonical isomorphism
$$
H^\dot(\M_g[l]) \cong H^\dot(\Gamma_g[l]).
$$
This holds with rational coefficients for all $l$, and arbitrary
coefficients whenever $\Gamma_g[l]$ is torsion free, which holds
whenever $Sp_g(\Z)[l]$ is torsion free --- $l \ge 3$.

We know from (\ref{van_h1}) that
$$
H^1(\Gamma_g[l]) \cong H^1(\M_g[l]) = 0
$$
when $g\ge 3$. By standard arguments (cf.\ \cite[\S 5]{hain:normal}),
this implies that
$$
c_1 : \Pic \M_g[l]\otimes\Q \to H^2(\M_g[l])
$$
is injective, and therefore that $\Pic \M_g[l]$ is finitely
generated when $g\ge 3$.

The stable cohomology of an arithmetic group depends only on
the ambient real algebraic group \cite{borel:triv}. Based on
this, one might expect that the natural map
$$
H^k(\G_g) \to H^k(\G_g[l])
$$
is an isomorphism for all $l\ge 0$, once the genus $g$ is sufficiently
large compared to the degree $k$. It follows from Johnson's work
that this is true when $k=1$ (cf.\ \cite{hain:normal}), but the only
evidence for it when $k>1$ is Harer's computation of the second homology
of the spin mapping class groups \cite{harer:spin}, and Foisy's theorem
from which Harer's computation now follows:

\begin{theorem}[Foisy \cite{foisy}]
For all $g\ge 3$, the natural map $H^2(\G_g) \to H^2(\G_g[2])$ is
an isomorphism. Consequently, $\Pic \M_g[2]$ is finitely generated
of rank 1.
\end{theorem}

\begin{question}
Is $\Pic \M_g[l]$ rank 1 for all $g\ge 3$ and all $l\ge 1$?
\end{question}

This would be the case if we knew that the $Sp_g(\Z)$ action
on $H^2(T_g)$ extended to an algebraic action of $Sp_g$,
for we could then invoke Borel's computation of the
stable cohomology of arithmetic groups \cite{borel:triv}.

\subsection{Normal Functions}
\label{subsec:normal} Each rational representation $V$ of
$Sp_g$ gives rise to an orbifold local system $\V$ over $\M_g[l]$. Such
a local system underlies an admissible variation of Hodge structure. First,
if $V$ is irreducible, then $V$ underlies a variation of Hodge
structure unique up to Tate twist (\cite[(9.1)]{hain:normal}).
Every polarized $\Q$ variation of Hodge structure whose monodromy
representation comes from a rational representation of $Sp_g$
has the property that each of its
isotypical components is an admissible variation of Hodge structure
of the form $A_\lambda\otimes \V(\lambda)$, where $A_\lambda$ is a
Hodge structure and $\V(\lambda)$ is a variation of Hodge structure
corresponding to the $Sp_g$ module with highest weight $\lambda$ ---
cf.\ \cite[(9.2)]{hain:normal}.

For a Hodge structure $V$ of weight $-1$ one
defines
the corresponding {\it intermediate jacobian} $JV$ by
$$
JV = V_\C/(F^0 V + V_\Z).
$$
Its interest comes from the fact that it parametrizes the extensions
of $\Z$ by $V$ in the MHS category:
if $E$ is an extension of the $\Z$ (with its trivial Hodge structure of
weight zero) by $V$, then choose an integral lift $e\in E$
of $1$ and consider the image of $e$ in
$$
E_\C/(F^0 E + V_\Z)\cong V_\C/(F^0 V + V_\Z).
$$
This is independent of the lift and yields a complete invariant of the
extension. There is an inverse construction that makes $JV$ support
a variation of mixed Hodge structure $\E$ that is universal as an
extension of the trivial Hodge structure $\Z$ by the constant
Hodge structure $V$:
$$
0 \to \V _{JV}\to \E \to \Z _{JV}\to 0
$$
(see \cite{carlson}). This immediately generalizes to a relative setting:
if $\V$ is an admissible variation of $\Z$ Hodge structure
of weight $-1$ over a smooth variety $X$, then we have a
corresponding bundle $\pi :\J\V\to X$ of intermediate jacobians over $X$
supporting a universal extension
$$
0 \to \pi ^*\V \to \E \to \pi ^*\Z _X\to 0.
$$
A section $\sigma$ of $\J\V$ over $X$ determines an extension of
Hodge structures:
$$
0 \to \V \to \sigma ^*\E \to \Z _X\to 0.
$$
A {\it normal function} is a section of $\J\V$ such that the corresponding
extension $\E$ is an admissible variation of mixed Hodge structure. The normal
functions arising from algebraic cycles are normal functions in this
sense  --- cf.\ \cite[\S6]{hain:normal}.

We briefly recall Griffiths' construction of a normal function
associated to a family of homologically trivial algebraic
cycles. First we consider the case where the base is a point.
Suppose that $X$ is a smooth projective variety. A homologically trivial
algebraic $d$-cycle $Z$ in $X$ canonically determines an extension of
$\Z$ by $H_{2d+1}(X;\Z(-d))$ by pulling back the exact sequence
$$
0 \to H_{2d+1}(X;\Z(-d)) \to H_{2d+1}(X,|Z|;\Z(-d)) \to H_{2d}(|Z|;\Z(-d))
\to \cdots
$$
of MHSs along the inclusion
$$
\Z \to H_{2d}(|Z|;\Z(-d))
$$
that takes 1 to the class of $Z$. So an integral lift of $1$ is
given by an integral singular $2d+1$ chain $W$ in $X$ whose boundary is $Z$.
Integration identifies $JH_{2d+1}(X;\Z (-d))$ with the {\it Griffiths
intermediate
Jacobian}
$$
J_d(X):=\Hom_\C (F^dH^{2d+1}(X);\C (-d))/H^{2d+1}(X;\Z (-d)),
$$
and under this  isomorphism
the extension class in question is just given by integration over $W$.

Families of homologically trivial cycles give rise to normal functions:
Suppose that $\X \to T$ is a family of smooth projective varieties
over a smooth base $T$ and that $\cZ$ is an algebraic cycle in
$\X$ which is proper over $T$ of relative dimension $d$. Then the local
system whose fiber over $t\in T$ is $H_{2d+1}(X_t;\Z(-d))$ naturally
underlies a variation of Hodge structure $\V$ over $T$ of weight $-1$
so that we can form the $d$th {\it relative intermediate jacobian}
$\J_d (\X /T)\to T$,
whose fiber over $t\in T$ is $J_d(X_t)$. The family of
cycles $\cZ$ defines a section of this bundle which is a normal function.

\begin{theorem}[Hain \cite{hain:normal}]\label{norm_classn}
Suppose that $\V$ is an admissible variation of Hodge structure
of weight $-1$ over $\M_g[l]$ whose monodromy representation factors
through a rational representation of $Sp_g$. If $g \ge 3$,
then the space of normal functions associated to $\V$ is finitely
generated of rank equal to the number of copies of the variation
$\V(\lambda_3)$ of weight $-1$ that occur in $\V$.
\end{theorem}

The theorem implies that, up to torsion and multiples, there is only
one normal function over $\M_g$ associated to a variation of Hodge
structure whose monodromy factors through a rational representation of
$Sp_g$. So what is the generator of these normal functions?

To answer this question, recall that if $C$ is a smooth projective curve of
genus $g$ and
$x\in C$, we have the Abel-Jacobi morphism
$$
C \to \Jac C, \quad y\mapsto (y)-(x).
$$
Denote the image 1-cycle in $\Jac C$ by $C_x$
and the cycle $i_\ast C_x$ by $C_x^-$, where $i : \Jac C \to \Jac C$
takes $u$ to $-u$. The cycle $C_x - C_x^-$ is homologous to zero, and
therefore defines a point $\nu ^1(C,x)$ in $J_1(\Jac C)$.
Pontrjagin product with the class of $C$ induces a homomorphism
$$
A:\Jac C \to J_1(\Jac C).
$$
We call the cokernel of $A$ the {\it primitive first intermediate Jacobian}
$J_1^{\pr}(\Jac C)$ of $\Jac C$. The family of such primitive intermediate
jacobians over $\M_g$ is the unique one (up to isogeny) associated to
the variation of Hodge structure of weight $-1$ over $\M_g$ whose associated
$\G_g$ module is $V(\lambda_3)$. It is not difficult to show that
$$
\nu ^1(C,x) - \nu ^1(C,y) = 2A(x-y).
$$

It follows that the image of $\nu ^1(C,x)$ in $J_1^{\pr}(\Jac C)$ is
independent of $x$. This is the value of the normal function associated with
$C-C^-$ over $[C]$.

We can do better and realize half of this generator by a generalized
normal function as follows.
Let $A$ be a principally polarized abelian variety of dimension $g\ge 3$.
The polarization determines a distinguished element $\omega$ of $H_2(A;\Z)$.
If $Z$ and $Z'$ are two piecewise smooth cycles representing
$\omega$, then their difference is the 
boundary of a piecewise smooth $3$-chain $W$ on $A$. 
Represent the dual of $H_3(A;\R)$ by  translation invariant $3$-forms on $A$.
Then integrating these forms over $A$ determines an element 
of $H_3(A;\R )$. Another choice of $W$ gives a  class
that differs from this one by an element of  $H_3(A;\Z )$, and so we have
a well-defined element $[Z-Z']$ of $H_3(A;\R/\Z )$. Notice that the latter torus
is naturally identified with the first intermediate jacobian $J_1(A)$ of $A$.
We declare $Z$ and $Z'$ to be equivalent if $[Z-Z']=0$ and denote the space of  
piecewise smooth cycles representing $\omega$ modulo this equivalence relation
by $D(A)$. This is clearly a torsor of $J_1(A)$ and so it has a natural
complex structure. In view of its connection with Deligne cohomology, we call
it the {\it Deligne torsor} of $A$.
This torsor contains naturally a subtorsor $D(A)[2]$ of the $2$-torsion in $J_1(A)$,
$J_1(A)[2]\cong H_3(A;\frac{1}{2}\Z/\Z )$: Let $a=(a_1,a_{-1},\dots ,a_g,a_{-g})$ be a symplectic basis of $H_1(A;\Z )$. Each basis element $a_i$ is uniquely represented by a homomorphism $\alpha_i:S^1\to A$ and so $\omega$ is represented by the $2$-cycle $\sum_{i=1}^g \alpha_i\times\alpha_{-i}$. This cycle defines an element $z(a)\in D(A)$. It is easily verified that $z(a)$ only depends on the mod two reduction of $a$ and that if $a$ runs over all symplectic bases, $z(a)$ runs over an entire orbit $D(A)[2]$ of $J_1(A)[2]$. (So $J_1(A)[2]\backslash D(A)$ has a canonical point which identifies it with $J_1(A)$.) The group $Sp(H_1(A;\Z))$ acts on $D(A)$ as an affine transformation group in a way that is easily made explicit. The lifts of these 
transformations to a universal covering of $D(A)$ form a group of affine symplectic transformations. It is an extension of $Sp(H_1(A;\Z))$ by $H_3(A;\Z )$ which splits if we enlarge the extension to $H_3(A;\frac{1}{2}\Z )$.

The Pontrjagin product with $\omega$ defines a homomorphism $A\to J_1(A)$ which gives rise to corresponding primitive notions: 
the {\it primitive Deligne torsor} $D^{\pr}(A):=A\backslash D(A)$ is a torsor
of the primitive intermediate Jacobian $J^\pr _1(A):=A\backslash J_1(A)$.
We have corresponding universal Deligne torsors over $\A_g$ which we denote 
$\cD _g\to\A_g$ and $\cD _g^{\pr}\to\A_g$. By the above argument, these torsors become trivial on the Galois cover of $\A_g$ representing principally polarized abelian varieties with a level $2$ structure. The torsors themselves are nontrivial, for it can be shown that the orbifold fundamental groups of these torsors are nonsplit extensions of the integral symplectic group of genus $g$

For $C$ a nonsingular projective curve of genus $g\ge 3$ and $x\in C$,
the Abel-Jacobi morphism $C \to \Jac C$ defined by $y\mapsto (y)-(x)$
defines a cycle
in the homology class of the natural polarization of $\Jac C$ and
so we get an element $[(C,x)]$ of $D(\Jac C)$. Its image in
$D^{\pr}(\Jac C)$ is independent of $x$ and so can be denoted by $[C]$.
Universally this produces holomorphic lifts of the period map: 
$$
\nu_g^1: \M_g^1\to \cD_g\text{  and  } \nu_g: \M_g\to \cD_g^{\pr}.
$$
We call $\nu_g$ the {\it fundamental normal function} on $\M_g$.

\subsection{Picard group of the generic curve with a level $l$
structure} \label{subsec:genericpicard}
The classification of normal functions (\ref{norm_classn}) implies that there
are no sections of $\Pic^0$ of infinite order defined over $\M_g[l]$ when $g\ge 3$.
This,  combined with Sipe's computation (\ref{sipe}) of the monodromy of 
roots of the canonical bundle allows one
to determine the Picard group of the generic point of $\M_g[l]$.
The case $l=1$ was the subject of the Franchetta Conjecture which
was deduced from Harer's computation of $\Gamma_g$ by Beauville
(unpublished) and by Arbarello and Cornalba \cite{arb_cor:pic}.

\begin{theorem}[Hain \cite{hain:derham}]
The Picard group of the generic curve of genus $g\ge 3$ with a level
$l$ structure is of rank 1, has torsion subgroup isomorphic to the
$l$ torsion points $H_1(\Jac S_g;\Z/l)$, and, modulo torsion, is
generated by the canonical bundle if $l$ is odd, and a theta
characteristic if $l$ is even.
\end{theorem}

\section{Relative Malcev Completion}
\label{sec:malcev}

Fundamental groups of smooth algebraic varieties are quite special
as we know from the work of Morgan \cite{morgan} and others.
The least trivial restrictions on these groups come from Hodge theory
and Galois theory. Since $\G_g$ is the (orbifold) fundamental group
of $\M_g$, a smooth orbifold, Hodge theory and Galois theory should
have something interesting to say about its structure. To put a MHS on a group
one needs to linearize it. One way to do this is to replace the
group by some kind of algebraic envelope and put a MHS on the
coordinate ring of this (pro)algebraic group.

In this section we introduce these linearizations and use them
to establish a relation between the fundamental normal function and a
remarkable central extension that is hidden in a
quotient of the mapping class group. Here the impact of mixed Hodge theory
is not yet felt, but we are setting the stage for
Section~\ref{sec:hodgemap} where it is omnipresent.

\subsection{Classical Malcev completion}
\label{subsec:malcev} Suppose that $\pi$ is a
finitely generated group. The classical Malcev (or unipotent) completion of
$\pi$
consists of a prounipotent group $\U(\pi)$ (over $\Q$) and a homomorphism
$\pi \to \U(\pi)$. It is characterized by the following universal mapping
property: if $U$ is a unipotent group, and $\phi : \pi \to U$ is a
homomorphism, there is a unique homomorphism of prounipotent groups
$\U(\pi) \to U$ through which $\phi$ factors. There are several well known
constructions of the unipotent completion, which can be found in
\cite{hain:comp}, for example. Each (pro)unipotent group $U$ is isomorphic
to its Lie algebra $\u$, a (pro)nilpotent Lie algebra via the exponential
map. Thus, to give the Malcev group $\U(\pi)$ associated to $\pi$ it suffices
to give its associated pronilpotent Lie algebra $\u(\pi)$. This Lie algebra is
called the {\it Malcev Lie algebra associated to $\pi$}. It comes with a
natural descending filtration whose $k$th term $\u ^{(k)}(\pi)$
is the closed ideal of $\u (\pi )$ generated by its $k$-fold commutators
$(k=1,2,\dots )$ and it is complete with respect to this filtration.
We will refer to this filtration as the {\it Malcev filtration}.

When $\pi$ is the fundamental group $\pi_1(X,x)$ of a smooth complex
algebraic variety, $\u(\pi)$ has a
canonical MHS which was first constructed by Morgan \cite{morgan}.
If $X$ is also complete, or more generally, when $H_1(X)$ has a pure Hodge
structure of weight $-1$, then the weight filtration is the Malcev
filtration:
$$
W_{-k}\u (\pi _1(X,x))=\u ^{(k)}(\pi _1(X,x)).
$$
Alternatively, this MHS determines and is determined by a MHS
on the coordinate ring $\O(\U(\pi))$ of the associated Malcev group.

We shall denote the Malcev completion of $\pi_g^n=\pi_1(S_g^n,x_0)$ by 
$\p_g^n$. 

\subsection{Relative Malcev completion}
\label{subsec:relmalcev} The Malcev completion
of a group $\pi$ is trivial when $H_1(\pi)$ vanishes, for
then $\pi$ has no non-trivial unipotent quotients. Since the
first homology of $\G_g$ vanishes for all $g$, its Malcev
completion will be trivial. Deligne has defined the notion
of Malcev completion of a group $\pi$ relative to a Zariski dense
homomorphism $\rho:\pi \to S$, where $S$ is a reductive algebraic
group defined over a base field $F$ (that we assume to be of
characteristic zero).

The {\it Malcev completion of $\pi$ relative to $\rho : \pi \to S$}
is a a proalgebraic $F$-group $\cG(\pi,\rho)$, which is an extension
$$
1 \to \U \to \cG(\pi,\rho) \to S \to 1
$$
of $S$ by a prounipotent group, together with a lift
$\rhotilde : \pi \to \cG(\pi,\rho)$ of $\rho$.\footnote{In many cases
the completion of $\pi$ over an algebraic closure $\bar F$ of $F$ is the
set of $\bar F$ points of the completion of $\pi$ over $F$. This is the
case for the mapping class groups when $g\ge 3$, but we do not know
whether this is true in general, except when $S$ is trivial.} It is
characterized by the following universal mapping property: if
$G$ is an $F$-group which is an extension of $S$ by a unipotent
group $U$, and if $\phi : \pi \to G$ is a homomorphism, then there
is a unique homomorphism $\cG(\pi,\rho) \to G$ through which $\phi$ factors:
$$
\phi : \pi \stackrel{\rhotilde}{\to} \cG(\pi,\rho) \to G.
$$
Since $S$ is reductive, we should think of $\U$ as the prounipotent radical
of $\cG (\pi ,\rho )$. One can show, for instance, that $\U$ has a Levi
supplement so that $\cG (\pi ,\rho )$ is a semidirect product of $S$ and
$\U$. The Lie algebra $\g (\pi ,\rho )$ of $\cG (\pi ,\rho )$ also comes
with a Malcev filtration with respect to which it is complete:
$\g (\pi ,\rho )^{(0)}=\g (\pi ,\rho )$, and for $k\ge 1$,
$\g (\pi ,\rho )^{(k)}$ is the closed ideal generated by $k$-fold
commutators in the Lie algebra of $\U$.

We will often write $\cG(\pi)$ instead of $\cG(\pi,\rho)$ when the
representation
$\rho$ is clear from the context. We shall denote the completion of
the (orbifold) fundamental group of a pointed orbifold $(X,x)$ with
respect to a Zariski dense reductive representation $\rho:\pi_1(X,x) \to S$
by $\cG(X,x;\rho)$, or simply $\cG(X,x)$ when $\rho$ is clear from the
context.

When $S$ is trivial, we recover the classical Malcev completion.
The universal property of the Malcev completion of $\ker \rho$ yields a
natural homomorphism of proalgebraic $F$-groups
$\U (\ker \rho)(F)\to \cG (\pi)$. In general, it is neither
surjective nor injective as the following two examples show.

\begin{example}
The fundamental group of the
symplectic Lie group $Sp_g(\R)$ is infinite cyclic and hence so is its
universal cover $\widehat{Sp}_g(\R)\to Sp_g(\R )$. This universal cover
is not an algebraic group (which follows for instance from the fact that
the complexification of $Sp_g(\R )$,  $Sp_g(\C )$, is simply connected).
The preimage $\widehat{Sp}_g(\Z)$ of $Sp_g(\Z)$ in this covering contains
the universal central extension of $Sp_g(\Z)$ by $\Z$. Now take for $\pi$
this central extension and for $\rho$ its natural homomorphism to
$Sp _g(\C )$. The corresponding relative Malcev completion is then reduced
to $Sp _g(\C )$ itself, so that the homomorphism from $\U (\Z)(\C)$ (which
is just the abelian group $\C$) to $\cG (\widehat{Sp}_g(\Z ))$ is
trivial. We will see that this example is realized inside a quotient of
the mapping class
group.
\end{example}

\begin{example}
In this example, $\ker \rho$ is trivial, but $\U(\pi)$ is not.
The basic fact we need (see \cite[(10.3)]{hain:derham}) is that there is
always a natural $S$ equivariant isomorphism
$$
H_1(\u(\pi)) \cong
\prod_{\alpha \in \Check{S}} H_1(\pi;V_\alpha)\otimes V_\alpha^\ast,
$$
where $V_\alpha$ denotes a representation with highest weight $\alpha$.
For $\pi$ we take $\G$, a finite index subgroup of $SL_2(\Z)$, for
$S$ we take $SL_2(\Q)$, and for $\rho$ we take the natural inclusion.
Denote the $n$th power of the fundamental representation of $SL_2$ by
$S^nV$.
For all such $\G$, there is an infinite number of integers $n\ge 0$ such
that $H^1(\G;S^nV)$ is non-trivial.\footnote{This is easily seen when
$\G$ is free, for example. In general it is related to the theory of
modular forms.} It follows that $\U(\G)$
has an infinite dimensional $H_1$, even though $\ker \rho$ is trivial.
\end{example}

This example suggests the following problem:

\begin{problem}
Investigate the relationship between the theory of modular forms
associated to a finite index subgroup $\G$ of $SL_2(\Z)$ and the
completion of $\G$ relative to the inclusion $\G \hookrightarrow
SL_2(\Q)$.
\end{problem}

\subsection{The relative Malcev completion of $\G_g$}
\label{malcevofgamma} The natural
homomorphism $\rho:\G_{g,r}^n \to Sp_g$ has Zariski dense image.
Denote the completion of $\G_{g,r}^n$ relative to $\rho$ by
$\cG_{g,r}^n$, its prounipotent radical by $\U_{g,r}^n$ and their Lie algebras
by $\g_{g,r}^n$ and $\u_{g,r}^n$.

The following theorem indicates the presence
of essentially one copy of the universal central extension of $Sp _g(\Z)$ in
quotients of each mapping class group of genus $g$ when $g\ge 3$.

\begin{theorem}[Hain \cite{hain:comp}] When $g\ge 2$, the homomorphism
\begin{equation}\label{nat_map}
\U(T_{g,r}^n) \to \U_{g,r}^n
\end{equation}
is surjective. When $g \ge 3$, its kernel
is a central subgroup isomorphic to the additive group.
\end{theorem}

This phenomenon is intimately related to the cycle $C-C^-$ and
its normal function as we shall now explain.

\subsection{The central extension}
\label{subsec:cent_extn} The existence of the central extension
has both a group theoretic and a geometric explanation. It is also
related to the Casson invariant through the work of Morita
\cite{morita:casson,morita:cocycles}. We begin with the group
theoretic one.

The group analogue of the Malcev filtration for the Torelli group $T_g$ is
the most rapidly descending central series of $T_g$ with
torsion free quotients:
$$
T_g = T_g^{(1)}\supset T_g^{(2)}\supset T_g^{(3)} \supset \cdots
$$
Note that $T_g^{(1)}/T_g^{(2)}$ is the maximal torsion free abelian
quotient of $T_g$, which is
$$
V(\lambda _3)_g\Z := \wedge^3 H_1(S_g;\Z)/\left([S_g]\times H_1(S_g;\Z )\right)
$$
by Johnson's Theorem (\ref{johnson}). The group $\G_g/T_g^{(3)}$
can be written as an extension
\begin{equation}
\label{ext1}
 1 \to T_g^{(2)}/T_g^{(3)} \to \G_g/T_g^{(3)}\to \G_g/T_g^{(2)} \to 1.
\end{equation}
It turns out that this sequence contains a multiple of the universal
central extension of $Sp_g(\Z )$ by $\Z$.

Since $V_g(\lambda _3)$ is a rational representation of
$Sp_g$, and since the surjection
$$
\wedge^2 V_g(\lambda _3)_\Z \to T_g^{(2)}/T_g^{(3)}
$$
induced by the commutator
is $Sp_g(\Z)$ equivariant, it follows that $T_g^{(2)}/T_g^{(3)}\otimes \Q$
is also a rational representation of $Sp_g$. Because $V_g(\lambda _3)$ is
an irreducible symplectic representation,
there is exactly one copy of the trivial representation in
$\wedge^2 V_g(\lambda _3)$. This copy of the trivial representation
survives in $T_g^{(2)}/T_g^{(3)}\otimes \Q$ \cite{hain:comp} so that there
is an $Sp_g(\Z)$ equivariant projection $T_g^{(2)}/T_g^{(3)}\to \Z$.
Pushing the extension (\ref{ext1}) out along this map gives an extension
\begin{equation}
\label{ext2}
0 \to \Z \to E \to \G_g/T_g^{(2)}\to 1
\end{equation}
Note that $E$ is a quotient of $\G_g$. We will manufacture a multiple
of the universal central extension of $Sp_g(\Z)$ from this group that
turns out to be the obstruction to the map $\U(T_g) \to \U_g$
being injective. (Full details can be found in \cite{hain:comp}.)

The group $\G_g/T_g^{(2)}$ can be written as an extension
\begin{equation}
\label{ext3}
 0 \to V_g(\lambda _3)_\Z \to \G_g/T_g^{(2)}\to Sp_g(\Z) \to 1.
\end{equation}
Morita \cite{morita:conj} showed that this
extension is {\it semisplit}, that is, if we replace $V_g(\lambda _3)_\Z$ by
$\frac{1}{2}V_g(\lambda _3)_\Z$, it splits.
(This can also be seen using the normal function of $C-C^-$.)

\begin{theorem}[Morita \cite{morita:cocycles}, Hain \cite{hain:comp}]
\label{nosplit}
The extension of $Sp_g(\Z)$ by $\Z$ obtained by pulling back the extension
(\ref{ext2}) along a semisplitting of (\ref{ext3}) contains
the universal central extension of $Sp_g(\Z)$.
\end{theorem}

The geometric picture uses the fundamental normal function $\nu_g$. 
The lifted period maps $\nu_g$ and $\nu_g^1$ to the Deligne torsors 
are easily seen to extend over the partial completions $\widetilde{\M}_g$ resp.\ $\widetilde{\M}^1_g$:
$$ 
\tilde\nu_g: \widetilde{\M}_g\to \cD_g^{\pr}, \quad 
\tilde\nu_g^1: \widetilde{\M}_g^1\to \cD_g.
$$ 
The orbifold fundamental group of $\cD_g^{\pr}$ resp.\ $\cD_{g,1}$
is an extension of $Sp_g(\Z )$ by $V_g(\lambda _3)_\Z$ resp.\ $\wedge ^3V_{g,\Z}$ 
as both the base and fiber are
Eilenberg MacLane spaces with these groups as orbifold fundamental groups. But by \ref{cor:pioftildem} the orbifold fundamental group of
$\widetilde{\M}_g$, resp.\ $\widetilde{\M}^1_g$, also has such a structure. 
Indeed: 

\begin{theorem}\label{fundgroup} 
For $g\ge 3$, the normal functions $\tilde\nu_g: \widetilde{\M}_g\to \cD_g^{\pr}$ and $\tilde\nu_g^1: \widetilde{\M}_g^1\to \cD_g$
induce an isomorphism on orbifold fundamental groups. (The former can be identified with $\G_g/T_g\!{}^{(2)}$ and the latter with
$\G_g^1/(T_g^1)^{(2)}$.) 
\end{theorem}

{}From this theorem we recover the fact that (\ref{ext3}) is semisplit, not split. But we get more, since it should also lead to a description of that extension.

The extension (\ref{ext2}) can also be realized geometrically.

\begin{proposition}[Hain \cite{hain:comp}]
There is a canonical (locally homogeneous) line bundle $\B_g$ over
the bundle $\cD_{g,1}^{\pr}\to\A_g$
that realizes the central extension (\ref{ext2}) via the isomorphism of the previous proposition
as an extension of orbifold fundamental groups. 
In particular, both $\tilde\nu^\ast \B_g$ and
$\nu^\ast \B_g$ have nonzero rational first Chern class.
The bundle $\nu^\ast \B_g$ is canonically metrized and
its square is isomorphic (as a metrized line bundle) to the metrized line bundle
associated to the archimedean height of the cycle $C-C^-$.
\end{proposition}

\section{Hodge Theory of the Mapping Class Group}\label{sec:hodgemap}

One reason that mixed Hodge theory is so powerful is that the MHS category is
abelian. In many situations this turns out to have
topological implications for algebraic varieties that are difficult, if
not impossible, to obtain directly. A somewhat
related (but less exploited) property is that a MHS
is canonically split over $\C$. This implies that the weight filtration
(which often has a topological interpretation) splits in a way that is
compatible
with all the algebraic structure naturally present. So, for many purposes,
there is no loss of information regarding this algebraic structure if we pass
to the corresponding weight graded object. For example, the Malcev
filtration on the Malcev Lie algebra of a smooth projective variety is minus
the weight filtration, and it therefore splits over $\C$ in a natural way.
This splitting is natural in the sense that it respects the Lie algebra
structure
and is preserved under all base point preserving morphisms. But if
we vary the complex structure on $X$ or the base point $x$, then
the splitting will, in general, vary  with it.

A basic example is the Malcev Lie algebra $\p_g^1$ of $\pi_g^1=\pi
_1(S_g^1,x_0)$.
The group $\pi _g^1$ is free on $2g$ generators and it is a classical
fact that the graded of $\p_g^1$ with respect to the
Malcev filtration is just the free Lie algebra generated by $V_g$. If
$S_g$ is given a conformal structure, then $V_g$ has a pure Hodge
structure of weight $-1$ and the weight filtration of $\p_g^1$
is minus the Malcev filtration. The splitting allows us to identify
$\p_g^1\otimes\C$ with the completion of $\Lie (V_g)\otimes\C$.

We shall come back to this example in Section~\ref{subsec:relation}. But
for now we will focus on the relative Malcev completions introduced in the
previous section.

\subsection{Hodge theory of $\cG_{g,r}^n$}
\label{subsec:hodge} A choice of a conformal structure on $S_g$
and nonzero tangent vectors at $x_{n+1},\dots ,x_{n+r}$, 
determines a point $x_o$ of the moduli space
$\M_{g,r}^n$. We can thus identify $\G_{g,r}^n$ with the
orbifold fundamental group of $(\M_{g,r}^n,x_o)$.
This induces an isomorphism of  $\cG_{g,r}^n$ with $\cG(\M_{g,r}^n,x_o)$,
the completion of $\pi_1(\M_{g,r}^n,x_o)$ with respect to the standard
symplectic representation. We shall write $\cG_{g,r}^n(x_o)$ for
$\cG(\M_{g,r}^n,x_o)$ and denote its prounipotent radical by $\U_{g,r}^n(x_o)$.
There is a general Hodge de~Rham theory of relative Malcev completion \cite{hain:derham}. Applying it to $(\M_{g,r}^n,x_o)$, one obtains the following
result:

\begin{theorem}[Hain \cite{hain:torelli}]\label{mhs}
For each choice of a base point $x_o$ of $\M_{g,r}^n$, there is a
canonical MHS on the coordinate ring $\O(\cG_{g,r}^n(x_o))$ which
is compatible with its Hopf algebra structure. Consequently, the Lie
algebra $\g_{g,r}^n(x_o)$ of $\cG_{g,r}^n(x_o)$ and
the Lie algebra $\u_{g,r}^n(x_o)$ of its prounipotent radical
both have a natural MHS.
\end{theorem}

Denote the Malcev Lie algebra of the subgroup of $\pi_1(\M_{g,r}^n,x_o)$
corresponding to the Torelli group $T_{g,r}^n$ by $\t_{g,r}^n(x_o)$.
The normal function of $C-C^-$ can be used to lift the MHS
from $\u_{g,r}^n(x_o)$ to $\t_{g,r}^n(x_o)$.

\begin{theorem}[Hain \cite{hain:torelli}]\label{cent_extn}
For each $g\ge 3$ and for each choice of a base point $x_o$ of
$\M_{g,r}^n$, there is a canonical MHS on $\t_{g,r}^n(x_o)$ which
is compatible with its bracket. Moreover, the canonical central
extension
$$
0 \to \Q(1) \to \t_{g,r}^n(x_o) \to \u_{g,r}^n(x_o) \to 0
$$
is an extension of MHSs, and the weight filtration equals the
Malcev filtration.
\end{theorem}

\subsection{A presentation of $\t_g$}
\label{subsec:torelli_presentn}
We denote the Malcev Lie algebra of $T_{g,r}^n$ by $\t_{g,r}^n$. The
existence of a MHS on $\t_{g,r}^n(x_o)$ implies that, after tensoring
with $\C$, there is a canonical isomorphism
$$
\t_{g,r}^n(x_o)\otimes \C \cong
\prod_m \Gr^W_{-m} \t_{g,r}^n(x_o)\otimes \C.
$$
Since the left hand side is (noncanonically) isomorphic to
$\t_{g,r}^n\otimes \C$, to give a presentation of $\t_{g,r}^n\otimes\C$,
it suffices to give a presentation of its associated graded.
It follows from Johnson's Theorem (\ref{johnson}) that each graded
quotient of the lower central series of $\t_g$ is a representation
of the algebraic group $Sp_g$. We will give a presentation of
$\Gr^W_\dot \t_g$ in the category of representations of $Sp_g$. Recall
that $\lambda_1,\dots, \lambda_g$ is a set of fundamental weights of
$Sp_g$. For a nonnegative integral linear combination of the fundamental
weights $\lambda =\sum _{i=1}^gn_i\lambda _i$ we denote by
$V_g(\lambda)$ the representation of $Sp_g$ with highest weight $\lambda$.

For all $g \ge 3$, the representation $\wedge^2 V_g(\lambda_3)$
contains a unique copy of $V_g(2\lambda_2) + V_g(0)$. Denote the
$Sp_g$ invariant complement of this by $R_g$. Since the quadratic
part of the free Lie algebra $\Lie (V_g)$ is $\wedge^2 V_g$, we can
view $R_g$ as being a subspace of the quadratic elements of
$\Lie (V_g(\lambda_3))$.

As mentioned earlier, it is unknown whether any $T_g$
is finitely presented when $g\ge 3$. But the following theorem says
that its de~Rham incarnation is:

\begin{theorem}[Hain \cite{hain:torelli}]\label{presentation}
For all $g \ge 3$, $\t_g$ is isomorphic to the completion of its
associated graded $\Gr^W_\dot\t_g$. When $g\ge 6$, this has presentation
$$
\Gr^W_\dot \t_g = \Lie (V_g(\lambda_3))/(R_g),
$$
where $R_g$ is the set of quadratic relations defined above.
When $3 \le g < 6$, the relations in $\Gr^W_\dot\t_g$ are generated by
the quadratic relations $R_g$, and possibly some cubic relations.
In particular, $\t_{g,r}^n$ is finitely presented whenever $g \ge 3$.
\end{theorem}

Note that this, combined with (\ref{cent_extn}) gives a presentation
of $\Gr^W_\dot\u_g$ when $g\ge 6$:

\begin{corollary}
For all $g \ge 3$, $\u_g$ is isomorphic to the completion of its
associated graded $\Gr^W_\dot\u_g$. When $g\ge 6$, this has quadratic
presentation
$$
\Gr^W_\dot \u_g = \Lie (V_g(\lambda_3))/(R_g + V_g(0)),
$$
where $R_g$ is the set of quadratic relations defined above and where
$V_g(0)$ is the unique copy of the trivial representation in 
$\wedge^2 V_g(\lambda_3)$.
When $3 \le g < 6$, the relations in $\Gr^W_\dot\u_g$ are generated by
the quadratic relations $R_g+V_g(0)$, and possibly some cubic relations.
In particular, $\u_{g,r}^n$ is finitely presented whenever $g \ge 3$.
\end{corollary}

The proof that the relations in the presentation of $\t_g$ are
generated by quadratic relations when $g\ge 6$ and quadratic
and cubic ones when $g \ge 3$ is not topological, but uses deep
Hodge theory and, surprisingly, intersection homology. The
key ingredients are a result of Kabanov, which we state below,
and M.~Saito's theory of Hodge modules.

We define the {\it Satake compactification} $\Mbar^{\rm sat}_g$ of
$\M_g$ as the closure of $\M_g$ inside the (Baily-Borel-)Satake
compactification of $\A_g$.

\begin{theorem}[Kabanov \cite{kabanov}]
For each irreducible representation $V$ of $Sp_g$, the
natural map
$$
IH^2(\Mbar^{\rm sat}_g;\V) \to H^2(\M_g;\V)
$$
is an isomorphism when $g\ge 6$. Here $\V$ denotes the generically
defined local system corresponding to $V$.
\end{theorem}

Such a local system $\V$ is, up to a Tate twist, canonically
a variation of Hodge structure. Saito's purity theorem then implies that
$H^2(\M_g;\V)$ is pure of weight $2+$ the weight of $\V$ when $g\ge 6$.
It is this purity result that forces $H^2(\t_g)$ to be of weight 2, and
implies that no higher order relations are needed.

\subsection{Understanding $\t_g$}\label{subsec:understanding}
Even though we have a presentation of $\t_g$, we still do not have a
good understanding of its graded quotients, either as vector spaces
or as $Sp_g$ modules.  There is an exact sequence
$$
 0 \to \p_g \to \t_g^1 \to \t_g \to 0
$$
of Lie algebras (recall that $\p_g$ stands for the Malcev Lie algebra
of $\pi _g=\pi_1(S,x_0)$).

It is the de~Rham incarnation of the exact sequence
of fundamental groups associated to the universal curve. Fix a conformal
structure on $(S,x_0)$.
Then this sequence is an exact sequence of MHSs.
Since $\Gr^W$ is an exact functor, and since
$\Gr^W_\dot \p_g$ is well understood, it suffices to
understand $\Gr^W_\dot \t_g^1$.

There is a natural representation
\begin{equation}\label{rep}
\t_g^1 \to \Der \p_g
\end{equation}
It is a morphism of MHS, and therefore determined by the graded
Lie algebra homomorphism
$$
\Gr^W_\dot\t_g^1 \to \Der \Gr^W_\dot \p_g.
$$
One can ask how close it is to being an isomorphism. Since this
map is induced by the natural homomorphism
\begin{equation}\label{can_hom}
\G_g^1 \to \varprojlim \Aut \C\pi_g/I^m,
\end{equation}
the homomorphism (\ref{rep}) factors through the projection
$\t_g^1 \to \u_g^1$,
and therefore cannot be injective. On the other hand, we have
the following (reformulated) result of Morita:

\begin{theorem}[Morita \cite{morita:trace}]
There is a natural Lie algebra surjection
$$
Tr_M : W_{-1}\Der \Gr^W_\dot \p_g \to \oplus_{k \ge 1} S^{2k+1}H_1(S)
$$
onto an abelian Lie algebra whose composition with (\ref{rep}) is
trivial. Here $S^m$ denotes the $m^{\text{th}}$ symmetric power.
\end{theorem}

One may then hope that the sequence
$$
0 \to \C \to \Gr^W_\dot\t_g^1 \to W_{-1}\Der \Gr^W_\dot \p_g \to
\oplus_{k \ge 1} S^{2k+1}H_1(S) \to 0
$$
is exact. However, there are further obstructions to exactness
at $W_{-1}\Der \Gr^W_\dot \p_g $ which were discovered by Nakamura
\cite{nakamura:obstn}.
They come from Galois theory and use the fact that $\M_g^1$ is defined
over $\Q$.\footnote{Actually, Nakamura proves his result for a
corresponding sequence for $\t_{g,1}$, but his obstructions most likely
appear in this case too.} On the other hand, one can ask:

\begin{question}
Is the map $\u_g^1 \to \Der \p_g$ injective? Equivalently,
is $\cG_g^1$ the Zariski closure of the image of the representation
(\ref{can_hom})?
\end{question}

A good understanding of $\t_g$ may help in understanding the
stable cohomology of $\G_g$ as we shall explain in the next
subsection.

\subsection{Torelli Lie algebras and the cohomology of $\G_g$}
\label{subsec:torelli&coho}
Each Malcev Lie algebra $\g$ can be viewed as a complete topological
Lie algebra. A basis for the neighbourhoods of 0 being the
terms $\g^{(k)}$ of the Malcev filtration. One can define the
continuous cohomology of such a $\g$ to be
$$
H^\dot(\g) := \varinjlim H^\dot(\g/\g^{(k)}).
$$
If $\g$ has a MHS, then so will $H^\dot(\g)$.
The continuous cohomology of $\t_{g,r}^n$, $\u_{g,r}^n$, etc. each has
an action of $Sp_g$. The general theory of relative Malcev completion
\cite{hain:derham} gives a canonical homomorphism
\begin{equation}\label{nat_homom}
H^\dot(\u_{g,r}^n)^{Sp_g} \to H^\dot(\M_{g,r}^n)
\end{equation}
One can ask how much of the cohomology of $\M_{g,r}^n$ is captured by
this map.

Fix a base point $x_o$ of $\M_{g,r}^n$. Then $\t_{g,r}^n$, etc.\ all
have compatible MHSs,
and these induce MHSs on their continuous cohomology
groups. These groups have the property that the weights on $H^k$ are
$\ge k$.

\begin{theorem}[Hain \cite{hain:torelli}]\label{morphism}
The map~(\ref{nat_homom}) is a morphism of MHS.
\end{theorem}

Since $H_1(T_g)$ is a quotient of $\u_g$, there is an induced map
\begin{equation}\label{ext_map}
H^\dot(H_1(T_g)) \to H^\dot(\u_g).
\end{equation}
This is also a morphism of MHS. The following result follows
directly from \cite[(9.2)]{hain:cycles}, the presentation
(\ref{presentation}) of $\t_g$, and the existence of the
MHS on $\u_g$.

\begin{proposition}
If $g\ge 3$, the map~(\ref{ext_map}) surjects onto the lowest weight
subring
$$
\oplus_{k\ge 0} W_kH^k(\u_g)
$$
of $H^\dot(\u_g)$, and the kernel is generated by the ideal generated
by the unique copy of $V_g(2\lambda_2)$ in $H^2(H_1(T_g))$.
\end{proposition}

Similar results hold when $r+n>0$ --- cf.\ \cite[\S 14.6]{hain:torelli}.

The following result of Kawazumi and Morita tells us that the
image of the lowest weight subring of $H^\dot(\u_g)^{Sp_g}$ contains
no new cohomology classes.

\begin{theorem}[Kawazumi-Morita \cite{kawazumi-morita}]
\label{kawazumi-morita}
The image of the natural map
$$
H^\dot(H_1(T_g))^{Sp_g} \to H^\dot(\M_g)
$$
is precisely the subring generated by the $\kappa_i$'s.
\end{theorem}

If we combine this with the previous two results and Pikaart's Purity
Theorem (\ref{pikaart_purity}), we obtain the following strengthening
of the theorem of Kawazumi and Morita (and obtained independently by Morita, 
building on our work):

\begin{theorem}\label{tautbound}
When $k\le g/2$, the image of $H^k(\u_g)^{Sp_g} \to H^k(\M_g)$
is the degree $k$ part of the subring generated by the $\kappa_i$'s.
\end{theorem}

To continue the discussion further, it seems useful to consider
cohomology with symplectic coefficients.

\subsection{Cohomology with symplectic coefficients}
The irreducible representations of $Sp_g$ are parametrized by Young
diagrams with $\le g$ rows (and no indexing of the boxes), in other
words, by nonincreasing sequences of nonnegative integers whose
terms with index $>g$ are zero.
So any such sequence $\alpha =(\alpha_1,\alpha_2,\dots)$ defines an
irreducible representation of $Sp_h$ for all $h\ge g$. We will
denote the representation of $Sp_g$ corresponding to $\alpha$
by $V_{g,\alpha}$, and the corresponding (orbifold) local system over
$\M_g$ by $\V_{g,\alpha}$.

A theorem of Ivanov \cite{ivanov} (that in fact pertains to more
general local systems) implies that, when $r\ge 1$, the group
$H^k(\G_{g,r}^n;V_{g,\alpha})$ is independent of $g$ once $g$ is
large enough. In the case at hand we have a more explicit result that we
state here for the undecorated case (a case that Ivanov actually excludes).

\begin{theorem}[Looijenga \cite{looijenga}]
Let $\alpha =(\alpha_1,\alpha_2,\dots)$ be a nonincreasing sequence of
nonnegative integers  that is eventually zero, and let $c_1,c_2,\dots$
be weighted variables with
$\deg (c_i)=2i$. Put $|\alpha |:=\sum _{i\ge 1} \alpha_i$.
Then there exists a finitely generated, evenly graded
$\Q [c_1,\dots ,c_{|\alpha |}]$-module $A^{\dot}_{\alpha}$ (that can be
described explicitly) and a
graded homomorphism of $H^\dot (\G _{\infty})$ modules
$$
A^\dot _{\alpha}[-|\alpha |]\otimes H^\dot(\G _{\infty})
\to H^\dot (\G_g;V_{g,\alpha })
$$
that is an  isomorphism in degree $\le cg-|\alpha |$. It is also
a MHS morphism if we take
$A^{2k}_{\alpha}$ to be pure of type $(k,k)$. In particular, we have
$$
A^\dot_{\alpha}[-|\alpha |]\otimes H^\dot(\G _{\infty})\cong
H^\dot(\G_\infty;V_\alpha )
$$
both as MHSs and as graded $H^\dot(\G _{\infty})$
modules. So, by (\ref{pikaart_purity}), $H^k(\G_\infty;V_\alpha)$ is
pure of weight $k+|\alpha|$.
\end{theorem}

It is useful to try to understand all cohomology groups with
symplectic coefficients at the same time. To do this we take a leaf
out of the physicist's book and consider the ``generating function''
\begin{equation}\label{gen_fn}
\oplus_\alpha H^\dot(\G_g;V_\alpha^\ast)\otimes V_\alpha
\end{equation}
where $\alpha$ ranges over all partitions with $\le g$ rows,
and $\ast$ denotes dual.
This is actually a graded commutative ring as the Peter-Weyl Theorem
implies that the coordinate ring $\O_g$ of $Sp_g$ is
$$
\O_g = \oplus_\alpha \left(\End V_\alpha\right)^\ast \cong
\oplus_\alpha V_\alpha^\ast\otimes V_\alpha.
$$
The mapping class group acts on $\O_g$ by composing the
right translation action of $Sp_g$ on $\O_g$ with the canonical
representation $\G_g \to Sp_g$. The corresponding cohomology group
$H^\dot(\G_g;\O)$ is then the ``generating function''~(\ref{gen_fn}).
Note that $\O_g$ is a variation of Hodge structure of weight 0, so
the group $H^k(\G_g;\O_g)$ is stably of weight $k$ by the above theorem.

There is a canonical algebra homomorphism
$$
H^\dot(\u_g) \to H^\dot(\G_g;\O_g)
$$
whose existence follows from the de~Rham theory of relative completion
suggested by Deligne --- cf.\ \cite{hain:derham}. The map~(\ref{nat_homom})
of the previous subsection is just its invariant part.
This map is a MHS morphism for each choice of complex
structure on $S$. The $\alpha$ isotypical
part of both sides stabilizes as $g$ increases. It is natural
to ask:

\begin{question}
Is this map stably an isomorphism?
\end{question}

This has been verified by Hain and Kabanov (unpublished) in degrees
$\le 2$ for all weights, and in degree 3 and weight 3.
If the answer is yes, or even if one has surjectivity, then it will
follow from the theorem of Kawazumi and Morita (\ref{kawazumi-morita})
that the stable cohomology of $\M_g$ is generated by the $\kappa_i$'s.
A consequence of injectivity
and Pikaart's Purity Theorem would be that for each $k$, $H^k(\u_g)$ is
pure of weight $k$ once the genus is sufficiently large. This
is equivalent to the answer to the following question being affirmative.

\begin{question}
Are $H^\dot(\u_g)$ and $U\Gr^W_\dot \u_g$ stably Koszul dual?
\end{question}

Note that $U\Gr^W_\dot \u_g$ and the lowest weight subalgebra of
$H^\dot(\u_g)$ have dual quadratic presentations.

\section{Algebras Related to the Cohomology of Moduli Spaces of Curves}
\label{sec:algebras}

The ribbon graph description is the root of a number of ways of
constructing (co)homology classes on moduli spaces of curves from certain
algebraic structures. These constructions have in common that they actually
produce cellular (co)chains on $\MM _g^n$, and so they
are recipes that assign numbers to `oriented' ribbon graphs. The typical
construction, due to Kontsevich \cite{kontsevich:feynman}, goes like this:
assume that we are given a complex vector space $V$,
a symmetric tensor $p\in V\otimes V$, and linear forms
$T_k:V^{\otimes k}\to \C$ that are cyclically invariant. If $\G$ is a
ribbon graph, then the decomposition of $X(G)$  into $\sigma _1$ and
$\sigma _0$-orbits gives isomorphisms
$$
\otimes _{s\in X_1(\G)}V^{\otimes\ori (s)}\cong V^{\otimes X(G)}\cong
\otimes _{v\in X_0(\G)}V^{\otimes\out (v)},
$$
where $\ori (s)$ stands for the two-element set of orientations of the edge
$s$ and $\out (v)$ for the set of oriented edges that have $v$ as initial
vertex. Now
$p^{\otimes X_1(\G)}$ defines a vector of the lefthand side
and a tensor product of certain $T_k$'s defines a linear form on the
righthand side. Evaluation of the linear form on the vector gives a number,
which is clearly an invariant of the ribbon graph.
Since this invariant does not depend on an orientation on the set of
edges of $G$, it cannot be used directly to define a cochain on the
combinatorial moduli spaces. To this end we need some sign rules so
that, for instance, the displayed isomorphisms acquire a sign.
The tensors $p$ and $T_k$ are sometimes referred to as the {\it propagator}
and the {\it interactions}, respectively,
to remind us of their physical origin.

If $p$ is nondegenerate, then we may use it to identify $V$ with its dual.
In this case $T_k$ defines a linear map $V^{\otimes (k-1)}\to V$.
The properties one needs to impose on propagator and interactions in order
that the above recipe produce a cocycle on $\MM _g^n$ is that they define a
$\Z /2$ graded $A_{\infty}$ algebra with inner product. A similar recipe
assigns cycles on $\MM_g^n$ to certain $\Z /2$ graded differential algebras.
The cocycles can be evaluated on the cycles and this, in principle, gives a
method of showing that some of the classes thus obtained are nonzero.
We shall not be more precise,
but instead refer to \cite{kontsevich:feynman} or \cite{seminar}
for an overview. A simple example is to take $V=\C$, $p:=1\otimes1$ and
$T_k(z^k)$ arbitrary for $k\ge 3$ odd, and zero otherwise. Kontsevich
asserts that the classes thus obtained are all tautological.

\subsection{Outer space}\label{subsec:outer}
In Section~\ref{sec:ribbon} we encountered a beautiful combinatorial
model for a virtual classifying space of the mapping class group $\G _g^n$.
There is a similar, but simpler, combinatorial model that does the same job
for the outer automorphism group of a free group.

We fix an integer $r\ge 2$ and consider connected graphs $G$ with first
Betti number equal to $r$ and where each vertex has degree
$\ge 3$. Let us call these graphs {\it $r$-circular graphs}.
The maximal number of edges (resp.\ vertices) such a graph can
have is $3r-3$ (resp.\ $2r-2$). These bounds are realized by all
trivalent graphs of this type. Notice that an $r$-circular graph $G$
has fundamental group isomorphic to the free group on $r$ generators, $F_r$.
We say that $G$ is {\it marked} if we are given an isomorphism
$\phi :F_r\to\pi _1(G,\text{base point})$ up
to inner automorphism. The group $\Out (F_r)$ permutes these markings
simply transitively. There is an obvious notion of isomorphism for marked
$r$-circular graphs. We shall denote the collection of isomorphism classes
by
$\cG _r$.

Let $(G,[\phi ])$ represent an element of $\cG _r$. The metrics on $G$ that
give $G$ total length $1$ are parameterized by the interior of a simplex
$\Delta (G)$. We fit these simplices together in a way analogous to the
ribbon graph case: if $s$ is an edge of $G$ that is not a loop, then
collapsing it defines another element $(G/s,[\phi ]/s)$ of $\cG _r$. We
may then identify $\Delta (G/s)$ with a face of $\Delta (G)$. After we have
made these identifications we end up with a simplicial complex $\OOhat _r$.
The union of the interiors of the simplices $\Delta (G)$ (indexed by $\cG _r$)
will be denoted by $\OO _r$; it is the complement of a closed subcomplex
of $\OOhat _r$. This construction is due to Culler-Vogtmann
\cite{culler}. We call $\OO _r$  the {\it outer space
of order $r$} for reasons that will become apparent in a moment.
Observe that $\OOhat _r$ comes with a simplicial action of
$\Out (F_r)$. We denote the quotient of $\OOhat _r$ (resp.\ $\OO _r$)
by $\Out (F_r)$ by $\GGhat_r$ (resp.\ $\GG _r$). It is easy to see that
$\GGhat _r$ is a finite orbicomplex. The open subset $\GG _r$ is the moduli
space of metrized $r$-circular graphs. It has a spine of dimension $2r-3$.

\begin{theorem}[Culler-Vogtmann \cite{culler}, Gersten]
The outer space of order $r$ is contractible and a subgroup of finite
index of
$\Out
(F_r)$ acts freely on it. Hence $\GG _r$ is a virtual classifying space for
$\Out (F_r)$ and $\Out (F_r)$ has virtual homological dimension $\le 2r-3$.
\end{theorem}

In contrast to the ribbon graph case, $\OO _r$ is not piecewise smooth. If we
choose $2g-1+n$ free generators for the fundamental group of our reference
surface $S_g^n$, then each ribbon graph without vertices of degree $\le 2$
determines an element of $\GG _{2g-1+n}$: simply forget the ribbon structure.
The ribbon data is finite and it is therefore not surprising that forgetting
the ribbons defines a finite map
$$
\widehat{f} : \cS _n\backslash\MMhat _g^n\to \GGhat _{2g-1+n}
$$
of orbicomplexes. Here $\cS _n$ stands for the symmetric group, which acts
in the obvious way on $\MMhat _g^n$. Following Strebel's theorem, the
preimage of
$\GG _{2g-1+n}$ can be identified with $\cS _n\backslash (\M _g^n\times
\inn \Delta ^{n-1})$. We denote the resulting map by
$$
f: \cS _n\backslash (\M _g^n\times \inn\Delta ^{n-1})\to \GG _{2g-1+n}.
$$
It induces the evident map
$$
f_*: H_k(\G _g^n)_{\cS _n}\to H_k(\Out F_{2g-1+n})
$$
on rational homology. It is unclear whether there is such an
interpretation for the induced map on cohomology with compact supports.
We remark that $\M _g^n\times \inn\Delta ^{n-1}$ is canonically oriented,
but that its $\cS _n$-orbit space is not (since transpositions reverse this
orientation). Poincar\'e duality therefore takes the form
$$
H_k(\M _g^n)_{\cS _n}\cong
H_k(\cS _n\backslash (\M _g^n\times \inn\Delta ^{n-1});\epsilon)\cong
H_c^{6g-7+3n-k}(\cS _n\backslash (\M _g^n\times \inn\Delta ^{n-1});
\epsilon ),
$$
where $\epsilon$ is the signum representation of $\cS _n$.
If $\delta$ denotes the (signum) character of $\Out (F_r)$ on
$\wedge ^rH_1(F_r)$, then the adjoint of $f_*$ is a map
$$
H_c^{6g-7+3n-k}(\GG _{2g-1+n};\delta )\to
H_c^{6g-7+3n-k}(\cS _n\backslash (\M _g^n\times \inn\Delta ^{n-1});
\epsilon)\cong H_k(\M _g^n)_{\cS _n}.
$$
So, when $m\ge 1$, we have maps
$$
H_c^{i+m-1}(\GG _{m+1};\delta )\to\oplus _{2g-2+n=m}H_{2m-i}
(\M _g^n)_{\cS _n}\to H_{2m-i}(\GG _{m+1}),\quad i=0,1,\dots .
$$
There is a remarkable interpretation of this sequence that we will discuss
next.

\subsection{Three Lie algebras}\label{subsec:threeLie}
We describe Kontsevich's three functors from the category of symplectic
vector spaces to the category of Lie algebras and their relation with
the cohomology of the moduli spaces $\M _g^n$. The basic
references are \cite{kontsevich:symp} and \cite{kontsevich:feynman}.

We start out with a finite dimensional $\Q$ vector space $V$ endowed with a
nondegenerate antisymmetric tensor $\omega _V\in V\otimes V$.
Let $\Ass (V)$ be the tensor algebra (i.e., the free associative algebra)
generated by $V$. We grade it by giving $V$ degree $-1$.  The Lie
subalgebra generated by $V$ is free and so we denote it by $\Lie (V)$.
It is well-known that $\Ass (V)$ may be identified with the universal
enveloping algebra of $\Lie (V)$. If we mod out $\Ass (V)$ by the two-sided
ideal generated by the degree $\le -2$ part of $\Lie (V)$,
we obtain the symmetric algebra $\Com (V)$ of $V$. Define
$\g _{\ass} (V)$ (resp.\ $\g _{\lie} (V)$) to be the Lie algebra of
derivations of $\Ass (V)$ (resp.\ $\Lie (V)$) of degree $\le 0$ that kill
$\omega _V$. Since each derivation of $\Lie (V)$ extends canonically to its
universal enveloping algebra, we have an inclusion
$\g _{\lie} (V)\subset\g _{\ass} (V)$. There is also a corresponding Lie
algebra $\g _{\com} (V)$ of derivations of degree $\le 0$ of $\Com (V)$
that kill $\omega _V$. Here we regard the latter as a two-form on the
affine space $\spec \Com (V)$. This Lie algebra is a quotient of $\g _{\ass}
(V)$. All three Lie algebras are graded and have as degree zero
summand the Lie algebra $\sp (V)$ of the group $\Symp (V)$ of symplectic
transformations of $V$. A simple verification shows that the degree $-1$ summands
have as $\sp (V)$ representations the following natural descriptions:
$$
\g_{\com}(V)_{-1}\cong S^3(V),\quad 
\g_{\ass}(V)_{-1}\cong S^3(V)\oplus \wedge ^3V,\quad 
\g_{\lie}(V)_{-1}\cong \wedge^3V.
$$
These Lie algebras are functorial with respect to symplectic injections
$(V,\omega _V)\hookrightarrow (W,\omega _W)$. Note that $\Symp (V)$ acts
trivially on this cohomology of the Lie algebra in question because
$\sp (V)\subset\g _*(V)$. This implies that $H^k(\g _*(V))$,
$\ast \in \{\lie,\ass,\com\}$, depends
only on $\dim V$. We form the inverse limit:
$$
H^k(\g _*):=\varprojlim_{V}H^k(\g _*(V)).
$$
The sum over $k$, $H^\dot(\g _*)$, has the structure of a connected graded bicommutative Hopf
algebra; the coproduct comes  from the direct sum operation on symplectic
vector spaces.
It is actually bigraded: apart from the cohomological grading
there is another coming from the grading of the Lie algebras. Notice that the
latter grading has all its degrees $\ge 0$. The primitive part
$H^{\dot}_{\pr}(\g _*)$ inherits this bigrading. Furthermore, the
natural maps
$$
H^{\dot}(\g _{\com})\to H^{\dot}(\g _{\ass})\to H^{\dot}(\g _{\lie})
$$
are homomorphisms of bigraded Hopf algebras. Consequently, we have induced
maps between the bigraded pieces of their primitive parts.

\begin{theorem}[Kontsevich \cite{kontsevich:symp}, \cite{kontsevich:feynman}]
\label{thm:lie}
For $\ast \in \{\lie,\ass,\com\}$ we have
$$
H^k_{\pr}(\g _*)_0=H^k_{\pr}(\sp_\infty )\cong
\begin{cases}
\Q &\text{for }k=3,7,11,\dots\, ;\\
0 &\text{otherwise}.
\end{cases}
$$
Furthermore,
$H^k_{\pr}(\g _*)_l=0$ when $l$ is odd and, when $m>0$, we have a natural
diagram
$$
\begin{CD}
H^k_{\pr}(\g _{\com})_{2m} @>>>
H^k_{\pr}(\g _{\ass})_{2m} @>>>
H^k_{\pr}(\g _{\lie})_{2m} \cr
@V{\cong }VV  @V{\cong }VV  @V{\cong }VV \cr
H^{k+m-1}_c(\GG _{m+1},\delta ) @>{f_c^*}>>
\oplus _{2g-2+n=m}H_{2m-k}(\M _g^n)_{\cS _n} @>{f_*}>>
H_{2m-k}(\GG_{m+1})\cr
\end{CD}
$$
which commutes up to sign and whose rows are complexes.
The maps in the top row are the natural maps and the bottom row is the
sequence defined in Section~\ref{subsec:outer}.
\end{theorem}

The proof is an intelligent application of classical invariant
theory. For each of the three Lie algebras one writes down the standard
complex. The subcomplex of invariants with respect to the symplectic
group is quasi-isomorphic to the full complex. Weyl's invariant theory
furnishes  a natural basis for this subcomplex. Kontsevich then
observes that this makes the subcomplex naturally isomorphic to a cellular
chain (or cochain) complex of one of the cell complexes $\GG _*$ and
$\MM _*^*$ whose (co)homology appears in the bottom row.

The diagram in this theorem suggests that  the
sequence of natural transformations $\Lie\to\Ass\to\Com$
is self dual in some sense. This
can actually be pinned down by looking at the corresponding operads:
Ginzburg and Kapranov \cite{ginz} observed that these operads have
``quadratic relations'' and they proved the self duality of the operad
sequence in a Koszul sense.

However, our main reason for displaying this diagram is that it pertains
to the cohomology of the moduli spaces of curves in two apparently
unrelated ways. The first one is evident. The $\cS _n$ coinvariants of
the homology of $\M _g^n$ features in the middle column, but the righthand
column has something to do with the cohomology of a
`linearization' of $\G _{\infty}$: we will see that  $\g _{\lie}$
is intimately related to the Lie algebra of the relative Malcev completions
discussed in Section~\ref{sec:malcev}. We explain this in the next subsection
after a giving a restatement of Kontsevich's Theorem.

In this restatement the Lie algebra cohomology of $\g_\ast (V)$ is replaced 
by the relative Lie algebra cohomology
of the pair $(\g_\ast (V),\k(V))$, where $\k(V)$ is a maximal compact Lie
subalgebra of $\g_\ast (V)_0$, and therefore of $\g_\ast (V)$ ($\k(V)$ is a
unitary Lie algebra of rank $\dim V/2$). As above, the Lie algebra cohomology
$H^k(\g_\ast (V),\k(V))$ depends only on the dimension of $V$ and stabilizes
once $\dim V$ is sufficiently large. We denote the inverse limit of these groups
by $H^\dot(\g_\ast,\k_\infty)$. Likewise, we denote the stable cohomology of the
pair $(\sp_g,\k_g)$ by $H^\dot(\sp_\infty,\k_\infty)$. By a theorem of Borel
\cite{borel:triv}, this is naturally isomorphic to the stable cohomology of $\A_g$
and is a polynomial algebra generated by classes $c_1, c_3, c_5, \dots$, where
$c_k$ has degree $2k$.

Combining Kontsevich's Theorem~\ref{thm:lie} with Borel's computation
and an elementary spectral sequence argument, we obtain the following
result. (Use the fact that $(\sp_g,\k_g)$
is both a sub and a quotient of $(\g_\ast(V),\k(V))$.)

\begin{corollary}\label{rel_lie}
We have
$$
H^k_{\pr}(\g_\ast,\k_\infty)_0=H^k_{\pr}(\sp_\infty,\k_\infty )\cong
\begin{cases}
\Q &\text{for }k=2,6,10,\dots\, ;\\
0 &\text{otherwise}.
\end{cases}
$$
Furthermore, when $m > 0$, the natural maps $H^\dot(\g_\ast,\k_\infty)_m
\to H^\dot(\g_\ast)_m$ are isomorphisms.
\end{corollary}

\subsection{Relation with the relative Malcev completion}
\label{subsec:relation}
We begin with an observation. For a symplectic vector space $V$ and
$\ast \in \{\lie,\ass,\com\}$, denote by $\g _*^{\flat}(V)$ the subalgebra of
$\g _*(V)$ generated by its summands of weight $0$ and $-1$.
Kontsevich's computation shows:

\begin{proposition} The graded cohomology groups
$H^k (\g _*^{\flat} (V),\k(V))_l$ stabilize and the sum of the stable terms is a bigraded
bicommutative Hopf algebra $H^\dot(\g_*^{\flat},\k_\infty)_{\dot}$. In addition,
the restriction map
$H^k_{\pr}(\g_*,\k_\infty)_l\to H^k_{\pr}(\g_*^{\flat}(V),\k_\lie(V))_l$
is an isomorphism when $l\le k$.
\end{proposition}

The case of interest here is that of lie  where 
$\g _\lie^{\flat}(V)_0=\sp (V)$ and $\g _\lie^{\flat}(V)_{-1}\cong \wedge ^3V$.
Denote by $z_m$ the element of $H^{2m}_{\pr}(\g_\lie )_{2m}$ that corresponds,
via Theorem~\ref{thm:lie}, to $1\in H_0(\GG _{m+1})$. The preceding proposition
yields:

\begin{corollary}\label{cor:weightcontrol}
We have a natural isomorphism of bigraded Hopf algebras
$$
\sum_{l\le k} H^k(\g_{\lie}^{\flat},\k_\infty)_l\cong
H^\dot(\sp_\infty,\k_\infty)[z_1,z_2,\dots ]
\cong \C[c_1,c_3,c_5,\dots,z_1,z_2,z_3,\dots]
$$
where each $z_i$ and $c_j$ is primitive.
\end{corollary}

The graded Lie algebra $\g^\flat_\lie(V)$ is the semi-direct product of
$\sp(V)$ and its elements of positive weight, which we shall denote by
$\u^\flat_\lie$. Consequently, there are natural inclusions
$$
H^\dot(\u^\flat_\lie(V))^{Sp} \hookrightarrow H^\dot(\g^\flat_\lie(V),\k(V))
\text{ and } H^\dot(\sp(V),\k(V)) \hookrightarrow
 H^\dot(\g^\flat_\lie(V),\k(V)).
$$
Together these induce an algebra homomorphism
$$
H^\dot(\sp(V),\k(V))\otimes H^\dot(\u^\flat_\lie(V))^{Sp}
\to H^\dot(\g^\flat_\lie(V),\k(V))
$$
which is compatible with stabilization.

\begin{proposition}\label{stab_prod}
Upon stabilization, these maps induce an isomorphism
$$
H^\dot(\sp_\infty,\k_\infty)\otimes H^\dot(\u^\flat_\lie)^{Sp}
\to H^\dot(\g^\flat_\lie,\k_\infty).
$$
\end{proposition}

Next we relate the graded Lie algebra $\g^\flat_\lie$ to the filtered
Lie algebra $\g_{g,1}$ of the relative Malcev completion $\cG _{g,1}$ of
$\G_{g,1}\to Sp_g$. 
Recall from Section~\ref{sec:groups} that $\pi_g^1$ is freely generated by $2g$
generators
named $\alpha _{\pm},\dots ,\alpha_{\pm g}$ so that 
the commutator $\beta:=(\alpha_1,\alpha_{-1})\cdots (\alpha_g,\alpha_{-g})$
represents a simple loop around $x_1$. Using Latin letters for the 
logarithms of the images of elements of $\pi _g^1$ in its Malcev completion,
we find that 
$$
b\equiv [a_1,a_{-1}]+\cdots +[a_g,a_{-g}]\mod{(\p _g^1)^{(3)}}.
$$
So the image of $b$ in $\Gr^2\p_g^1\cong \wedge^2 V_g$ is the symplectic
form $\w_S$. 

The obvious homomorphism $\G_{g,1} \to \Aut (\pi_g^1)$ induces a Lie algebra
homomorphism
\begin{equation}\label{map}
\g _{g,1} \to \Der \p_g^1.
\end{equation}
whose image we denote by $\gbar _{g,1}$.\footnote{It is possible that
this map is injective so that $\g _{g,1}\cong\gbar _{g,1}$. Note that it
is not surjective --- see \cite{morita:trace} and \cite{nakamura:obstn},
and Section~\ref{subsec:understanding}.}

Notice that $\gbar_{g,1}$ is contained in the subalgebra $\Der(\p_{g,1},b)$
consisting of those derivations that kill $b$.  Since (\ref{map}) is (Malcev)
filtration preserving, it induces Lie
algebra homomorphisms
$$
\Gr^\dot \g_{g,1}\to \Gr^\dot\gbar_{g,1}\to \Der^\dot(\Gr\p_g^1,\w_S).
$$
Notice that the last term is just $\g _{\lie}(V_g)$.

In view of (\ref{stab_prod}), to construct an algebra homomorphism
$$
H^\dot(\g^\flat_\lie,\k_\infty) \to H^\dot(\Gamma_\infty)
$$
it suffices to construct an algebra homomorphism
$H^\dot(\u^\flat_\lie)^{Sp_g} \to H^\dot(\Gamma_\infty)$.

At this stage we need Hodge theory. Choose a conformal structure
on $S_g$. Then, by~(\ref{mhs}), there are 
natural MHSs on $\g_{g,1}$ and $\p_{g,1}$ whose weight filtrations are the Malcev filtrations and such that (\ref{map}) is a MHS morphism. Hence 
the image $\gbar_{g,1}$ has a natural MHS. 
Since $\Gr^\dot \g_{g,1}$ is generated by
its summands in degree $0$ and $1$, the same is true for $\Gr^\dot\gbar_{g,1}$.
On the other hand, the summands in degree $0$ and $1$ of 
$\Gr^\dot\gbar_{g,1}$ are equal to the summands of weight $0$ and $-1$ of $\g
_{\lie}(V_g)$, and so the graded Lie algebra $\Gr^\dot\gbar_{g,1}$
may be identified with $\g_\lie^\flat (V_g)_{\dot}$ (except that the indexing 
of the summands differs by sign).

Denote the pronilpotent radical $W_{-1}\gbar_{g,1}$ of $\gbar_{g,1}$ by
$\ubar_{g,1}$. We know from Section~\ref{sec:hodgemap} that the homomorphisms
\begin{equation}\label{sequence}
H^\dot (\ubar_{g,1})^{Sp_g}\to H^\dot(\u_{g,1})^{Sp_g}\to H^\dot(\Gamma_{g,1})
\end{equation}
are morphisms of MHS. After weight grading these become bigraded algebra 
homomorphisms
$$
H^\dot(\u_\lie^\flat (V_g)_{\dot})^{Sp_g}\to H^\dot(\Gr^W_{\dot}\u_{g,1})^{Sp_g}
\to \Gr^W_\dot H^\dot(\G_{g,1}).
$$
The sequence (\ref{sequence}) stabilizes with $g$ to a sequence of Hopf
algebras in the MHS category. The corresponding weight graded sequence is
$$
H^\dot((\u_\lie^\flat )_{\dot})^{Sp_g}\to 
H^\dot(\Gr^W_{\dot}\u_{\infty,1})^{Sp_g}
\to \Gr^W_\dot H^\dot(\G_{\infty}).
$$
Each term in this sequence is a Hopf algebra and each map a Hopf algebra
homomorphism. But by Pikaart's Purity Theorem we know that the last term
is pure of weight $k$ in degree $k$, so that we can replace it by
$H^\dot(\G_{\infty})$ and obtain a map $H^\dot(\u^\flat_\lie)^{Sp_g}
\to H^\dot(\Gamma_\infty)$. We therefore have a Hopf algebra homomorphism
$$
H^\dot(\g^\flat_\lie,\k_\infty) \cong
H^\dot(\sp_\infty,\k_\infty)\otimes H^\dot(\u^\flat_\lie)^{Sp_g}
\to H^\dot(\Gamma_\infty).
$$

If we compose the  natural restriction map
$H^\dot(\g_\lie, \k_\infty)\to H^\dot(\g_\lie^\flat,\k_\infty)$ with the
above maps we get a homomorphism
$$
H^\dot(\g_\lie,\k_\infty)\to H^\dot(\G_{\infty}).
$$
Kontsevich asked (at the end of \cite{kontsevich:feynman}) about the meaning 
of this map.\footnote{Actually, Kontsevich asks this for with $H^\dot(\g_\lie)$
in place of the relative Lie algebra cohomology. However, there
does not seem to be a natural homomorphism to $H^\dot(\G_\infty)$ in this
case.} This can now be answered by invoking the theorem of Kawazumi and Morita 
(\ref{kawazumi-morita}), or rather a weaker form, which says that $z_i$ is
mapped to a nonzero multiple of $\kappa_i$. This result was obtained
with Kawazumi and Morita.

\begin{theorem}
There is a natural Hopf algebra homomorphism
$$
H^\dot(\g_\lie,\k_\infty)\to H^\dot(\G_{\infty}).
$$
The left hand side is a polynomial algebra generated by primitive 
elements $z_1, z_2,\dots$ and $c_1, c_3,  \dots $ where $z_i$ has degree
$2i$ and $c_j$ has degree $2j$. The image of this homomorphism is 
precisely the subalgebra generated by the $\kappa_i$s. The kernel is
generated by elements of the form
$$
c_{2k+1} - a_k z_{2k+1} - P_{2k+1}, \quad k \in \{0,1,2,\dots \}
$$
where $P_{2k+1}$ is a polynomial in the $z_i$ and $c_j$ with no linear
terms , and $a_k$ is a non zero rational number. 
\end{theorem}

Here we have used the fact, due to Mumford \cite{mumford}, that the image
of $c_{2k+1}$ in $H^\dot (\M_{g,1})$ is a polynomial in the odd $\kappa_i$'s.
The theorem indicates that no new stable classes are to be expected from
Hodge theory --- that is, a de~Rham version of the Mumford conjecture
holds. Recent work of Kawazumi and Morita attempts to explain the kernel 
of the homomorphism $H^{\dot}(\gbar_{\infty})\to H^{\dot}(\G_{\infty })$
in terms of secondary characteristic classes of surface bundles. The first 
element of the kernel is the difference $c_1 - 12\, z_1$. Its restriction
to the Torelli group can be interpreted as the Casson Invariant (cf.\ \cite{morita:casson}.)

\end{document}